\documentclass[11pt]{amsart}
\usepackage[T1]{fontenc}
\usepackage{lmodern}
\usepackage{mathtools,amssymb,esint,microtype}
\usepackage[a4paper,margin=25mm]{geometry}
\usepackage{etoolbox,needspace}
\makeatletter
\patchcmd{\@settitle}{\uppercasenonmath\@title}{}{}{
  \errmessage{Title case patch failed}}
\patchcmd{\@setauthors}{\MakeUppercase{\authors}}{\authors}{}{
  \errmessage{Author case patch failed}}
\patchcmd{\author@andify}{\unskip ,\penalty-2}{\unskip\penalty-2}{}{
  \errmessage{Author comma patch failed}}
\patchcmd{\maketitle}{\uppercasenonmath\shorttitle}{}{}{
  \errmessage{Running title case patch failed}}
\patchcmd{\maketitle}{\@nx\MakeUppercase{\the\toks@}}{\the\toks@}{}{
  \errmessage{Running author case patch failed}}
\makeatother

\usepackage[hidelinks,bookmarksnumbered]{hyperref}
\hypersetup{pdftitle={Gradient potential estimates for k-Hessian equations},pdfauthor={Huajie Liu, Shiguang Ma, Shengyang Zang},pdfsubject={Gradient potentials, gradient-growth equations, and sharpness for Hessian equations}}
\numberwithin{equation}{section}
\newtheorem{theorem}{Theorem}[section]
\newtheorem{lemma}[theorem]{Lemma}
\newtheorem{proposition}[theorem]{Proposition}
\newtheorem{corollary}[theorem]{Corollary}
\theoremstyle{definition}
\newtheorem{definition}[theorem]{Definition}
\theoremstyle{remark}
\newtheorem{remark}[theorem]{Remark}
\newcommand{\R}{\mathbb R}
\newcommand{\W}{\mathcal W}
\newcommand{\restr}{\mathbin{\vrule height 1.4ex depth -0.3ex width .07ex\vrule height .07ex depth 0ex width .7ex}}
\DeclareMathOperator{\Lip}{Lip}
\DeclareMathOperator{\osc}{osc}
\DeclareMathOperator{\tr}{tr}
\DeclareMathOperator{\supp}{supp}
\newcommand{\avint}{\fint}
\allowdisplaybreaks[2]
\title{Gradient potential estimates for $k$-Hessian equations}
\date{}
\author{Huajie Liu}
\address[Huajie Liu]{School of Mathematical Sciences, Peking University,
Beijing, 100871, China}
\email{huajieliu@math.pku.edu.cn}

\author{Shiguang Ma}
\address[Shiguang Ma]{School of Mathematical Sciences, Nankai University,
Tianjin, 300071, China}
\email{msgdyx8741@nankai.edu.cn}

\author{Shengyang Zang}
\address[Shengyang Zang]{School of Mathematical Sciences, Nankai University,
Tianjin, 300071, China}
\email{1120230032@mail.nankai.edu.cn}

\begin{document}
\begin{abstract}
For $2\le k<n$, we prove pointwise gradient bounds for proper weak
$k$-convex functions in terms of a local average and a Wolff
potential of the full Hessian measure $F_k[u]$.
The estimate yields local Lipschitz regularity for $F_k[u]=f$
with $f\in L^{n/k,1/k}_{\mathrm{loc}}$, and for
$F_k[u]=|Du|^p$ when $0<p<nk/(n-k)$.
Both the Lorentz index $1/k$ and the upper bound on $p$ are sharp.
\end{abstract}

\maketitle
\setcounter{tocdepth}{1}
\tableofcontents

\section{Introduction}\label{sec:introduction}

We study pointwise gradient estimates for $k$-Hessian equations with
measure data. For a smooth function $u$, the $k$-Hessian operator
$\sigma_k(D^2u)$ is the $k$-th elementary symmetric function of the
eigenvalues of $D^2u$. Trudinger--Wang \cite{TW1,TW} extended this
operator to a nonnegative Radon measure $F_k[u]$ for proper weak
$k$-convex functions and proved its weak continuity under local
$L^1$ convergence. This formulation includes Hessian measures with
singular parts.

Estimates for function values include H\"older bounds under density
or ball-mass assumptions \cite{Tru97,TW} and Wolff-potential estimates
due to Labutin \cite{Lab} and Phuc--Verbitsky \cite{PV}.
Interior gradient estimates for smooth admissible solutions were
obtained by Trudinger \cite{Tru97} and Chou--Wang \cite{CW}.
The gradient potential theory for quasilinear equations
\cite{DM10,KuM,KuM14,NP} suggests a corresponding question for
Hessian equations: how does the full Hessian measure control the
gradient of a weak $k$-convex function?

For $2\le k<n$, we bound the gradient of a proper weak
$k$-convex function $u$ at Lebesgue points of $Du$ by a local
average of $u$ and a Wolff potential of its full Hessian measure.
In a comparison argument over successive scales, ordinary interior
gradient estimates may introduce a fixed multiplier at each step.
We prove a Lipschitz stability estimate for $k$-harmonic functions
in which the preceding Lipschitz constant enters with coefficient
one. The comparison errors can then be summed and bounded by
the potential.

The mapping properties of potentials \cite{Cianchi} then yield
local Lipschitz regularity for $F_k[u]=f$ with
$f\in L^{n/k,1/k}_{\mathrm{loc}}$. We also prove local Lipschitz
regularity for $F_k[u]=|Du|^p$ when $0<p<nk/(n-k)$.
The second Lorentz index $1/k$ and the upper bound on $p$ are
sharp in the full weak class.

\subsection{Main results}

For a symmetric matrix $A$, let $\sigma_j(A)$ denote the $j$-th
elementary symmetric function of its eigenvalues, with $\sigma_0=1$,
and set
\[
\Gamma_k=\{A:\sigma_j(A)>0,\ 1\le j\le k\}.
\]
Let $\Omega\subset\R^n$ be open.

\begin{definition}\label{def:weakclass}
	\textup{(\cite{TW})}
	For $1\le k\le n$, a function $u:\Omega\to[-\infty,\infty)$ is
	\emph{proper weak $k$-convex} if
	\begin{enumerate}
		\item $u$ is upper semicontinuous;
		\item $u$ is not identically $-\infty$ on any connected component
		of $\Omega$;
		\item every $C^2$ function $\varphi$ touching $u$ from above at
		$x\in\Omega$ 
        satisfies
		$D^2\varphi(x)\in\overline{\Gamma_k}$.
	\end{enumerate}
	We denote this class by $\Phi_k(\Omega)$.
\end{definition}

For $u\in\Phi_k(\Omega)$, write $\mu=F_k[u]$ for its full
Hessian measure, following the notation in
\cite[Section~7, p.~903]{PV}. For smooth $u$,
\[
F_k[u]=\sigma_k(D^2u)\,dx.
\]
Throughout the paper, $F_k[u]=f$ means that the full Hessian
measure equals $f\,dx$. For $2\le k<n$, weak gradients are locally
integrable by \eqref{eq:sobolev}, and their pointwise values are
understood as Lebesgue values. We write
$\avint_E=|E|^{-1}\int_E$.

For a nonnegative Radon measure $\mu$, define
\begin{equation*}
	\W_k^R[\mu](x)=\int_0^R
	\left(\frac{\mu(B_t(x))}{t^{n-k}}\right)^{1/k}\frac{dt}{t},
	\quad 0<R<\infty.
\end{equation*}
In the usual notation, this is
$\mathbf W_{k/(k+1),\,k+1}^R[\mu]$.

\begin{theorem}\label{thm:main}
	Let $2\le k<n$, let $\Omega\subset\R^n$ be open, and let
	$u\in\Phi_k(\Omega)$. Set $\mu=F_k[u]$.
	If $B_R(x_0)\Subset\Omega$ and $x_0$ is a Lebesgue point of $Du$,
	then
	\begin{equation}\label{eq:main}
		|Du(x_0)|\le\frac CR\avint_{B_R(x_0)}|u|
		+C\W_k^R[\mu](x_0),
	\end{equation}
	where $C=C(n,k)$.
\end{theorem}

\begin{remark}
	The measure $F_k[u]$ may contain a singular part, and $u$ may be
	unbounded. Subtracting an affine function preserves both
	$k$-convexity and the Hessian measure. Thus, for
	$\ell(x)=a\cdot(x-x_0)+b$, applying \eqref{eq:main} to $u-\ell$
	gives
	\begin{equation}\label{eq:affine}
		|Du(x_0)-a|\le\frac CR\avint_{B_R(x_0)}|u-\ell|
		+C\W_k^R[\mu](x_0).
	\end{equation}
	Choosing $b=\avint_{B_R(x_0)}u$ and using Poincar\'e's inequality,
	we obtain
	\begin{equation}\label{eq:gradientmean}
		|Du(x_0)-a|\le C\avint_{B_R(x_0)}|Du-a|
		+C\W_k^R[\mu](x_0).
	\end{equation}
	For $k=1$, the gradient estimate follows from linear potential
	theory, and $\W_1^R$ is a truncated first-order Riesz potential,
	up to normalization. For $k=n$, convexity gives interior gradient
	bounds in terms of local averages.
\end{remark}

Applying H\"older's inequality to the potential in
Theorem~\ref{thm:main} gives the following local Lipschitz estimate.

\begin{corollary}\label{cor:lqlipschitz}
	Let $2\le k<n$, $u\in\Phi_k(\Omega)$, and
	$B_R(x_0)\Subset\Omega$. Suppose that $F_k[u]=f$ in $B_R(x_0)$,
	with $f\in L^q(B_R(x_0))$ for some $n/k<q\le\infty$. Set
	\[
	\alpha_q=1-\frac{n}{kq}>0\text{ and } \alpha_\infty=1.
	\]
	Then, for every affine function $\ell$ with $D\ell=a$,
	\begin{equation}\label{eq:lqlipschitz}
		\|Du-a\|_{L^\infty(B_{R/2}(x_0))}
		\le \frac{C}{R}\avint_{B_R(x_0)}|u-\ell|\,dx
		+\frac{C}{\alpha_q}R^{\alpha_q}
		\|f\|_{L^q(B_R(x_0))}^{1/k},
	\end{equation}
	where $C=C(n,k)$. For $q=\infty$, the second term on the
	right-hand side reduces to
	$CR\|f\|_{L^\infty(B_R(x_0))}^{1/k}$. In particular, if $F_k[u]=f$ throughout $\Omega$, then
	\begin{equation}\label{eq:lqcriterion}
		f\in L^q_{\mathrm{loc}}(\Omega),\quad n/k<q\le\infty
		\quad\Longrightarrow\quad
		u\in W^{1,\infty}_{\mathrm{loc}}(\Omega).
	\end{equation}
\end{corollary}

If $F_k[u]=f$ in $\Omega$, the critical Lorentz condition
$f\in L^{n/k,1/k}_{\mathrm{loc}}(\Omega)$ gives local Lipschitz
regularity by Proposition~\ref{data:critical-lipschitz}.
The example in Proposition~\ref{prop:criticaldensity} shows that
the second Lorentz exponent $1/k$ is sharp.
For $\max\{0,2-n/k\}<\beta<1$, the function $u(x)=|x|^\beta$
is bounded and $k$-convex in $B_1$, while
$|Du(x)|=\beta|x|^{\beta-1}\to\infty$ as $x\to0$.

We next use Theorem~\ref{thm:main} to obtain local Lipschitz
regularity for equations with gradient growth.
Related equations for smooth solutions were studied
by Mohammed--Porru \cite{MP26}.

\begin{theorem}\label{growth:lipschitz}
	Let $2\le k<n$, $0\le\Lambda<\infty$, and
	$0\le h\in L^q_{\mathrm{loc}}(\Omega)$ for some $n/k<q\le\infty$.
	Suppose that $u\in\Phi_k(\Omega)$ satisfies
	\begin{equation}\label{growth:equation}
		F_k[u]=f\quad\text{in }\Omega,\quad
		0\le f\le\Lambda|Du|^p+h\quad\text{almost everywhere in }\Omega,
	\end{equation}
	with
	\begin{equation}\label{growth:range}
		0<p<p_*:=\frac{nk}{n-k}.
	\end{equation}
	Then $u\in W^{1,\infty}_{\mathrm{loc}}(\Omega)$.
	
	In particular, every weak $k$-convex solution of
	$F_k[u]=|Du|^p$ in $\Omega$ is locally Lipschitz whenever
	$p$ satisfies \eqref{growth:range}.
\end{theorem}

The proof of Theorem~\ref{growth:lipschitz} is given in
Section~\ref{sec:gradient-growth}. The restriction $p<p_*$
is sharp among weak $k$-convex solutions:
Proposition~\ref{growth:critical-example} constructs a solution
of $F_k[u]=|Du|^{p_*}$ in $B_1$ whose gradient belongs to
$L^{p_*}(B_1)$ but is unbounded. The example is continuous
when $n<2k$ and unbounded when $n\ge2k$. In the latter range,
the endpoint case for locally bounded solutions falls outside
this construction.

Section~\ref{sec:potential-scale} also gives estimates for the
oscillation of $u$, including a log-Lipschitz modulus under a
uniform ball-mass bound $\mu(B_r(x))\le Mr^{n-k}$.
The gradient can be discontinuous: for $k\ge2$, the function
$u(x)=|x_1|$ satisfies $F_k[u]=0$.
Gradient continuity is available under stronger hypotheses.
Caffarelli \cite{Caf91} proved interior $C^{1,\alpha}$ regularity
for strictly convex Alexandrov solutions of the Monge--Amp\`ere
equation with density bounded between positive constants.
Chen--Wang--Wu \cite[Theorem~1.4]{CWW23} obtained
$C^{1,\alpha}_{\mathrm{loc}}$ regularity for convex weak
$k$-Hessian solutions in a bounded convex domain, with zero or
affine boundary data and $0<c_0\le f\le c_1$.

The \hyperref[sec:examples]{appendix} gives examples with vanishing
critical Morrey density and gradients outside
$\mathrm{BMO}_{\mathrm{loc}}$. For $k\ge3$, gradient BMO can also
fail for densities in $L^{n/k}$.

\medskip
\noindent\textbf{Note added.}
The pointwise gradient potential estimates in this paper and in
He~\cite{He26} were obtained independently. In the common range
$2\le k<n$, the Wolff-potential estimate in
\cite[Theorem~1.1]{He26} corresponds to \eqref{eq:gradientmean} with
$a=0$, up to fixed changes of the radii. Both proofs use a Lipschitz
stability estimate with coefficient one for $k$-harmonic functions.

\subsection{Proof strategy and organization}

The proof of Theorem~\ref{thm:main} begins with a two-point maximum
argument based on the negative-direction Newton tensor inequality
\cite[(3.5)]{Tru97}. For
$v\in\Phi_k(B_r)\cap C(\overline B_r)$ with $F_k[v]=0$ and
$L,\delta\ge0$, we prove
\[
|v(x)-v(y)|\le L|x-y|+\delta\quad(x,y\in B_r)
\quad\Longrightarrow\quad
\Lip(v;B_{r/2})\le L+C(n)\frac{\delta}{r}.
\]
An integral estimate gives the initial Lipschitz bound.
We then restore the Hessian measure on successive annuli and
apply this stability estimate to the comparison functions.
The coefficient of $L$ remains one, so the errors add over scales.
The zero-boundary estimate of Phuc--Verbitsky \cite{PV} controls
these errors when $n>2k$; dimensional extension handles $n\le2k$.
An averaging argument gives the estimate for weak solutions.
Remark~\ref{rem:trudinger-comparison} explains the role of the
stability estimate in this summation.

Section~\ref{sec:mainproof} contains the proof of the pointwise
estimate and Corollary~\ref{cor:lqlipschitz}.
Section~\ref{sec:gradient-growth} treats equations with gradient
growth, and Section~\ref{sec:potential-scale} develops the
regularity consequences. The appendix contains sharpness
constructions and examples of gradient BMO failure.

\section{The pointwise gradient estimate}\label{sec:mainproof}
In this section we prove Theorem~\ref{thm:main} and
Corollary~\ref{cor:lqlipschitz}. We begin with the basic properties
of weak $k$-convex functions and the potential estimates needed below.

\subsection{Preliminaries}\label{sec:preliminaries}

We write $B_r=B_r(0)$ and $\omega_m=|B_1^m|$.
Unless explicitly fixed, constants may vary from line to line,
with their dependence indicated in each statement.
Norms of weak gradients are understood in the essential sense.
Inequalities between measures are understood on all Borel sets.
Measures appearing in potentials are extended by zero outside
their domains when necessary.

The approximation and comparison arguments below use the following
properties of weak $k$-convex functions.

\begin{lemma}
	\label{lem:hessian-tools}\label{lem:superadd}
	\textup{(\cite{TW1,TW,TW3})}
	Let $2\le k<n$ and let $\Omega\subset\R^n$ be a bounded domain.
	\begin{enumerate}
		\renewcommand{\labelenumi}{\textup{(\roman{enumi})}}
		\item
		We have
		\begin{equation}\label{eq:sobolev}
			\Phi_k(\Omega)\subset W^{1,q}_{\mathrm{loc}}(\Omega),
			\quad 1\le q<\frac{nk}{n-k}.
		\end{equation}
		
		\item
		For $u\in\Phi_k(\Omega)$, its standard mollifications
		$u_\varepsilon$ are smooth and $k$-convex on each subdomain
		$\Omega'\Subset\Omega$ for sufficiently small $\varepsilon>0$,
		and $u_\varepsilon\to u$ in $L^1_{\mathrm{loc}}(\Omega)$.
		Moreover, if $u_m\in\Phi_k(\Omega)$ and
		$u_m\to u$ in $L^1_{\mathrm{loc}}(\Omega)$, then $u$ admits
		a unique representative in $\Phi_k(\Omega)$. For this representative,
		\begin{equation}\label{eq:weakcont}
			F_k[u_m]\to F_k[u]
			\quad\text{locally weakly as Radon measures}.
		\end{equation}
		
		\item
		If $a,b\in\Phi_k(\Omega)$, then $a+b\in\Phi_k(\Omega)$ and
		\begin{equation*}
			F_k[a+b]\ge F_k[a]+F_k[b].
		\end{equation*}
		
		\item
		For every ball $B\subset\R^n$ and
		$a,b\in\Phi_k(B)\cap C(\overline B)$,
		\begin{equation*}
			\left.
			\begin{aligned}
				F_k[a]&\ge F_k[b]\quad\text{in }B,\\
				a&\le b\quad\text{on }\partial B
			\end{aligned}
			\right\}
			\quad\Longrightarrow\quad a\le b\quad\text{in }B.
		\end{equation*}
	\end{enumerate}
\end{lemma}

The homogeneous comparison uses $k$-harmonic functions,
which we define in the admissible viscosity sense.

\begin{definition}\label{def:kharmonic}
	\textup{(\cite{TW})}
	For $1\le k\le n$, a function
	$u\in C(\Omega)\cap\Phi_k(\Omega)$ is called
	\emph{$k$-harmonic} in $\Omega$ if the following conditions
	hold at every $x\in\Omega$:
	\begin{enumerate}
		\item
		Every $C^2$, $k$-convex test function $\varphi$ 
        touching $u$ from below at $x$ satisfies
		\[
		\sigma_k(D^2\varphi(x))\le0.
		\]
		\item
		Every $C^2$ test function $\varphi$ 
		touching $u$ from above at $x$ satisfies
		\[
		\sigma_k(D^2\varphi(x))\ge0.
		\]
	\end{enumerate}
\end{definition}

For $u\in C(\Omega)\cap\Phi_k(\Omega)$, these conditions
are equivalent to $F_k[u]=0$ in $\Omega$; see \cite{TW}.
We construct the comparison functions by solving Dirichlet
problems on balls.

\begin{lemma}\label{lem:bounded-dirichlet}
	\textup{(\cite{Tru97,TW})}
	Let $2\le k<n$ and let $B\subset\R^n$ be a ball.
	Suppose that $\varphi\in\Phi_k(B)\cap C(\overline B)$ and
	$0\le f\in L^\infty(B)$. Then there is a unique
	$v\in\Phi_k(B)\cap C(\overline B)$ satisfying
	\begin{equation*}
		\begin{cases}
			F_k[v]=f &\text{in }B,\\
			v=\varphi &\text{on }\partial B.
		\end{cases}
	\end{equation*}
	If $f=0$, then $v$ is locally Lipschitz in $B$.
	If $f\in C^\infty(\overline B)$,
	$\varphi|_{\partial B}\in C^\infty(\partial B)$, and
	$\inf_B f>0$, then $v$ is a classical solution.
\end{lemma}

To control the difference between successive comparison functions,
we use a potential estimate for solutions with zero boundary data.

\begin{lemma}\label{lem:value}
	\textup{(\cite{Lab,PV})}
	Let $2\le k<n/2$ and $s>0$.
	Let $\eta$ be a nonnegative measure on $B_{2s}$ with bounded
	density and $\supp\eta\subset\overline{B_s}$.
	Suppose that $z\in\Phi_k(B_{2s})\cap C(\overline{B_{2s}})$
	solves
	\[
	\begin{cases}
		F_k[z]=\eta &\text{in }B_{2s},\\
		z=0 &\text{on }\partial B_{2s}.
	\end{cases}
	\]
	Then, for every $x\in B_{2s}$,
	\begin{equation}\label{eq:value}
		0\le -z(x)\le C(n,k)\int_0^{8s}
		\left(\frac{\eta(B_t(x))}{t^{n-2k}}\right)^{1/k}
		\frac{dt}{t}.
	\end{equation}
	Here $\eta$ is extended by zero outside $B_{2s}$.
\end{lemma}

In the passage to weak solutions, the weak convergence in
\eqref{eq:weakcont} will be combined with the following
$L^1$ convergence of truncated potentials.

\begin{lemma}\label{lem:potentialcontinuity}
	Let $2\le k<n$ and let $G\subset\R^n$ be open.
	Suppose that $\nu_m$ and $\nu$ are nonnegative Radon measures
	on $G$ and that $\nu_m\to\nu$ locally weakly in $G$.
	Let $E\subset G$ be bounded and measurable, and let $a>0$.
	Assume that there is a compact set $K\subset G$ such that
	$\overline{B_a(y)}\subset K$ for every $y\in E$.
	Then
	\begin{equation*}
		\W_k^a[\nu_m]\to\W_k^a[\nu]
		\quad\text{in }L^1(E).
	\end{equation*}
\end{lemma}

\begin{proof}
	Choose a compact neighbourhood of $K$ contained in $G$.
	Local weak convergence gives a common bound $M$ for the masses
	of $\nu_m$ and $\nu$ on this neighbourhood.
	H\"older's inequality and Fubini's theorem give, for
	$0<\epsilon<a$,
	\begin{equation*}
		\sup_m\int_E\W_k^\epsilon[\nu_m](x)\,dx
		\le C|E|^{1-1/k}M^{1/k}\epsilon,
	\end{equation*}
	and the same estimate holds for $\nu$.
	
	For each $x\in E$, weak convergence gives
	$\nu_m(B_t(x))\to\nu(B_t(x))$ for all but countably many
	$t\in[\epsilon,a]$, since only countably many spheres
	$\partial B_t(x)$ have positive $\nu$-measure.
	Dominated convergence on $E\times[\epsilon,a]$, followed by
	$\epsilon\to0$, proves the claim.
\end{proof}

We also need to control an average of the potential by its
value at the center, with a larger truncation radius.

\begin{lemma}\label{lem:averagepotential}
	Let $2\le k<n$, and let $\nu$ be a nonnegative Radon measure
	defined in a neighbourhood of $\overline{B_R(x_0)}$.
	For every $0<\rho<R/4$,
	\begin{equation*}
		\avint_{B_\rho(x_0)}\W_k^{R/2}[\nu](y)\,dy
		\le C(n,k)\W_k^R[\nu](x_0).
	\end{equation*}
\end{lemma}

\begin{proof}
	Split the potential integral at $t=\rho$.
	Concavity and Fubini's theorem give
	\[
	\avint_{B_\rho(x_0)}\W_k^\rho[\nu](y)\,dy
	\le C\left(
	\frac{\nu(B_{2\rho}(x_0))}{\rho^{n-k}}
	\right)^{1/k}.
	\]
	The right-hand side is bounded by a constant times the
	contribution of $[2\rho,4\rho]$ to $\W_k^R[\nu](x_0)$.
	
	For $\rho\le t\le R/2$ and $y\in B_\rho(x_0)$, we have
	$B_t(y)\subset B_{2t}(x_0)$. The change of variables
	$t\mapsto2t$ therefore bounds the remaining part by
	$C\W_k^R[\nu](x_0)$.
\end{proof}

\subsection{Lipschitz stability for \(k\)-harmonic functions}\label{sec:homogeneous}\label{subsec:negative-direction}
The homogeneous comparison rests on an algebraic perturbation in a
direction where the Hessian is nonpositive. We first record this
ingredient and then apply it in a two-point maximum argument.
The cone $\Gamma_k$ is convex and invariant under orthogonal conjugation and addition of nonnegative matrices. The function $\sigma_k^{1/k}$ is concave on this cone. If $A\in\Gamma_k$, its restriction to a hyperplane belongs to $\Gamma_{k-1}$. The Newton tensors
\[
 T_{j-1}(A)=D_A\sigma_j(A),\quad 1\le j\le k,
\]
are positive definite and satisfy
\begin{equation}\label{eq:trace}
 \tr T_{j-1}(A)=(n-j+1)\sigma_{j-1}(A).
\end{equation}
The corresponding non-strict statements hold on the closed cones. These algebraic properties are recalled in \cite[Section~2]{Wang}.

The negative-direction estimate \eqref{eq:negative} follows the
Newton tensor argument in \cite[Section~3]{Tru97}, applied also to
lower-order operators. We prove the form needed here and derive
the perturbation conclusion \eqref{eq:strictincrease}.

\begin{lemma}\label{lem:matrix}
Let $2\le k<n$, $A\in\Gamma_k$ and $|e|=1$. If $A[e,e]\le0$, then
\begin{equation}\label{eq:negative}
 T_{j-1}(A)[e,e]\ge\sigma_{j-1}(A),\quad 1\le j\le k.
\end{equation}
Consequently, if $b>0$ and $A[e,e]+nb<0$, then, with
\[
 V=b\bigl((n+1)e\otimes e-I\bigr),
\]
we have
\begin{equation}\label{eq:strictincrease}
 A+V\in\Gamma_k,\quad \sigma_k(A+V)>\sigma_k(A).
\end{equation}
\end{lemma}
\begin{proof}
Choose the first coordinate in the direction $e$ and write
\[
 A=\begin{pmatrix}a&\xi^T\\ \xi&B\end{pmatrix},\quad a\le0.
\]
Then $B\in\Gamma_{k-1}$ and $T_{j-1}(A)[e,e]=\sigma_{j-1}(B)$. For $j\ge3$,
expansion in principal minors gives
\[
 \sigma_{j-1}(A)=\sigma_{j-1}(B)
       +a\sigma_{j-2}(B)-\xi^TT_{j-3}(B)\xi
       \le\sigma_{j-1}(B).
\]
The cases $j=1,2$ follow from $T_0=I$ and $\tr A=a+\tr B$, respectively.

Set $A_t=A+tV$, $0\le t\le1$. Since $A_t[e,e]=A[e,e]+nbt<0$,
\eqref{eq:trace} and \eqref{eq:negative} imply, as long as $A_t\in\Gamma_k$,
\begin{equation}\label{eq:path}
 \frac{d}{dt}\sigma_j(A_t)
 =b\bigl((n+1)T_{j-1}(A_t)[e,e]-\tr T_{j-1}(A_t)\bigr)
 \ge jb\,\sigma_{j-1}(A_t)>0,\quad 1\le j\le k.
\end{equation}
At a first exit from $\Gamma_k$, all the $\sigma_j$ would still be at least
their strictly positive initial values, a contradiction. Thus the entire
segment lies in $\Gamma_k$, and integrating \eqref{eq:path} with $j=k$
proves \eqref{eq:strictincrease}.
\end{proof}

The interior Lipschitz estimate for $k$-harmonic functions follows
from \cite[Theorem~3.1]{Tru97} by approximation. For the iteration
below, we need the more precise estimate \eqref{eq:defectbound},
in which the coefficient of $L$ is one. This allows us to sum
the errors from successive comparisons, as in \eqref{eq:iteration}.

\begin{proposition}\label{prop:defect}
Let $2\le k<n$ and let $v\in\Phi_k(B_r)\cap C(\overline{B_r})$ be
$k$-harmonic. There is a constant $C_0=C_0(n)$ such that $v$ is Lipschitz in $B_{r/2}$; more precisely,
\[
 \Lip(v;B_{r/2})\le \frac{C_0}r\osc_{B_r}v.
\]
Moreover, whenever $L,\delta\ge0$ satisfy
\begin{equation}\label{eq:defect}
 |v(x)-v(y)|\le L|x-y|+\delta\quad(x,y\in B_r),
\end{equation}
one has the sharper estimate
\begin{equation}\label{eq:defectbound}
 \Lip(v;B_{r/2})\le L+C_0\frac{\delta}{r}.
\end{equation}
In particular, if $w$ is Lipschitz in $B_r$ and $0\le w-v\le\delta$, then
\begin{equation}\label{eq:LS}
 \Lip(v;B_{r/2})\le\Lip(w;B_r)+C_0\frac{\delta}{r}.
\end{equation}
\end{proposition}
The interior Lipschitz estimate follows from \eqref{eq:defectbound}
by taking $L=0$ and $\delta=\osc_{B_r}v$.
\begin{proof}
Suppose first that $v$ is smooth and admissible, with
$\sigma_k(D^2v)=c>0$ in $B_r$, so that $D^2v\in\Gamma_k$.
The case $\delta=0$ is immediate from \eqref{eq:defect}.
For $\delta>0$, fix $x_0\in B_{r/2}$ and set
\[
 K=32(n+1),\quad
 \phi(s)=Ls+K\delta\left(\frac{2s}{r}-\frac{s^2}{r^2}\right),
\]
for $0\le s\le r$.
Then
\begin{equation*}
  \phi'(s)=L+2K\delta\left(\frac{1} r-\frac{s}{r^2}\right),
 \quad \phi''(s)=-\frac{2K\delta}{r^2}
\end{equation*}
and
\begin{equation}\label{eq:phi-above-linear}
 \phi(s)\ge Ls
\end{equation}
On $\overline{B_{ r/2}(x_0)}\times\overline{B_{r/2}(x_0)}$ consider
\[
 \Psi(x,y)=v(x)-v(y)-\phi(|x-y|)
 -\frac{8\delta}{r^2}\bigl(|x-x_0|^2+|y-x_0|^2\bigr).
\]
By \eqref{eq:defect} and \eqref{eq:phi-above-linear},
\begin{equation}\label{eq:Psi-boundary}
 \Psi(x,y)\le-\delta<0
 \quad\text{if }x\text{ or }y\in\partial B_{r/2}(x_0),
\end{equation}
whereas
\begin{equation}\label{eq:Psi-diagonal}
 \Psi(x,x)=-\frac{16\delta}{r^2}|x-x_0|^2\le0.
\end{equation}
We claim that
\begin{equation}\label{eq:gradient-at-x0}
 |Dv(x_0)|\le\phi'(0)=L+\frac{2K\delta}{r}.
\end{equation}
Otherwise, with $e_0=Dv(x_0)/|Dv(x_0)|$, Taylor's formula gives
\[
 \Psi(x_0+te_0,x_0)
 =t\bigl(|Dv(x_0)|-\phi'(0)\bigr)+o(t)>0
\]
for sufficiently small $t>0$. Thus \eqref{eq:Psi-boundary} and
\eqref{eq:Psi-diagonal} place a positive maximum at an interior pair
$(\bar{x},\bar{y})$ with $\bar{x}\ne \bar{y}$.

At this pair set
\[
 e=\frac{\bar{x}-\bar{y}}{|\bar{x}-\bar{y}|},\quad A=D^2v(\bar{x}),\quad B=D^2v(\bar{y}),
\]
and let
\[
 H=\phi''e\otimes e+\frac{\phi'}{|\bar{x}-\bar{y}|}(I-e\otimes e),
\]
with the derivatives of $\phi$ evaluated at $|\bar{x}-\bar{y}|$.
The second-derivative test yields
\begin{equation}\label{eq:doublehessian}
 \begin{pmatrix}A&0\\0&-B\end{pmatrix}
 \le\begin{pmatrix}H&-H\\-H&H\end{pmatrix}
      +\frac{16\delta}{r^2}I_{2n}.
\end{equation}
Testing on $(e,0)$ gives the negative direction
\begin{equation}\label{eq:direction}
 A[e,e]\le\phi''+\frac{16\delta}{r^2}
 =\left(-2K+16\right)\frac{\delta}{r^2}.
\end{equation}
Next let $Q=I-2e\otimes e$ be the orthogonal reflection across $e^\perp$.
Since
\begin{equation*}
 z-Qz=2(z\cdot e)e, \forall z\in\R^n,
\end{equation*}
testing \eqref{eq:doublehessian} on $(z,Qz)$ removes the tangential part of
$H$ and gives
\[
 z^T(A-QBQ)z
 \le4\phi''(z\cdot e)^2+\frac{32\delta}{r^2}|z|^2.
\]
Equivalently,
\begin{equation}\label{eq:reflection}
 QBQ\ge A+\frac{8K\delta}{r^2}e\otimes e
                -\frac{32\delta}{r^2}I.
\end{equation}
Set $b=32\delta/r^2$ and $V=b((n+1)e\otimes e-I)$. By
\eqref{eq:direction} and the choice of $K$,
\[
 A[e,e]+nb\le-(32n+48)\frac{\delta}{r^2}<0.
\]
Lemma~\ref{lem:matrix} therefore implies
\begin{equation}\label{eq:matrix-strict-here}
 A+V\in\Gamma_k,\quad\sigma_k(A+V)>\sigma_k(A).
\end{equation}
Moreover, $2K\ge8(n+1)$, so \eqref{eq:reflection} gives
\begin{equation*}
 P:=QBQ-(A+V)\ge0.
\end{equation*}
The segment $C_t=A+V+tP$, $0\le t\le1$, remains in $\Gamma_k$, and
$\frac{d}{dt}\sigma_k(C_t)=\tr(T_{k-1}(C_t)P)\ge0$. Consequently,
\begin{equation*}
 \sigma_k(QBQ)\ge\sigma_k(A+V).
\end{equation*}
Orthogonal invariance and \eqref{eq:matrix-strict-here} now give
\[
 c=\sigma_k(B)=\sigma_k(QBQ)\ge\sigma_k(A+V)>\sigma_k(A)=c,
\]
a contradiction. Thus \eqref{eq:gradient-at-x0} holds for every
$x_0\in B_{r/2}$. Integrating along line segments gives
\begin{equation}\label{eq:smooth-defect-bound}
 \Lip(v;B_{r/2})\le L+C_0\frac{\delta}{r},
\end{equation}
with the fixed choice $C_0=2K$, independent of $c$.

Next, a sequence of smooth solutions is used to approximate the homogeneous weak solution. Choose $\varphi_m\in C^\infty(\partial B_r)$ with
$\epsilon_m:=\|\varphi_m-v\|_{L^\infty(\partial B_r)}\to0$, and constants
$c_m>0$ tending to zero. The classical Dirichlet theorem
\cite[Theorem~3.2]{TW} on the uniformly $(k-1)$-convex ball gives smooth
admissible solutions
\begin{equation*}
 \begin{cases}
 F_k[v_m]=c_m&\text{in }B_r,\\
 v_m=\varphi_m&\text{on }\partial B_r.
 \end{cases}
\end{equation*}

Define
\[
 \alpha_m=\left(\frac{c_m}{\binom nk}\right)^{1/k},
 \quad q_m(x)=\frac{\alpha_m}{2}(|x|^2-r^2).
\]
Then
\begin{equation*}
 \begin{cases}
 F_k[q_m]=c_m&\text{in }B_r,\\
 q_m=0&\text{on }\partial B_r,
 \end{cases}
\end{equation*}
and $q_m\le0$ in $B_r$.

By Lemma~\ref{lem:superadd}\textup{(iii)},
$F_k[v+q_m]\ge c_m\,dx=F_k[v_m]$, while
$v+q_m-\epsilon_m\le\varphi_m$ on $\partial B_r$. Comparison gives
\begin{equation*}
 v+q_m-\epsilon_m\le v_m\quad\text{in }B_r.
\end{equation*}
Conversely, $F_k[v_m]\ge F_k[v+\epsilon_m]=0$ and
$\varphi_m\le v+\epsilon_m$ on the boundary, so
\begin{equation*}
 v_m\le v+\epsilon_m\quad\text{in }B_r.
\end{equation*}
Thus
\begin{equation*}
 v+q_m-\epsilon_m\le v_m\le v+\epsilon_m,
\end{equation*}
and, since $\|q_m\|_\infty=\alpha_mr^2/2\to0$,
\begin{equation*}
 \|v_m-v\|_{L^\infty(B_r)}\to0.
\end{equation*}
The approximants satisfy
\[
 |v_m(x)-v_m(y)|\le L|x-y|+\delta+2 \|v_m-v\|_{L^\infty(B_r)}
.
\]
Applying \eqref{eq:smooth-defect-bound} and passing to the uniform limit
proves \eqref{eq:defectbound}.

Taking $L=0$ and $\delta=\osc_{B_r}v$ gives the ordinary interior estimate.
Finally, if $0\le g:=w-v\le\delta$, then
$|g(x)-g(y)|\le\delta$, so \eqref{eq:defect} holds with
$L=\Lip(w;B_r)$. This proves \eqref{eq:LS}.
\end{proof}

The following interpolation consequence of Proposition~\ref{prop:defect} supplies the initial term in the
proof of Theorem~\ref{thm:main}; see also
\cite[Corollary~3.4]{TW} and \cite[Lemma~4.1]{Wang}.

\begin{corollary}\label{cor:L1Lip}
Let $2\le k<n$. If $h\in\Phi_k(B_r)\cap C(\overline{B_r})$ is $k$-harmonic, then, for every $a\in\R$,
\begin{equation}\label{eq:L1Lip}
 \Lip(h;B_{r/4})\le\frac{C(n)}r\avint_{B_r}|h-a|.
\end{equation}
\end{corollary}
\begin{proof}
Replacing $h(x)$ by $h(rx)-a$, it suffices to prove
\begin{equation}\label{eq:L1Lip-unit}
 \Lip(h;B_{1/4})\le C(n)\avint_{B_1}|h|.
\end{equation}
With $C_1=2C_0$, Proposition~\ref{prop:defect} gives
\begin{equation}\label{eq:oscLip}
 \Lip(h;B_{t/2}(y))\le\frac{C_1}t\sup_{B_t(y)}|h|\text{ for all } \overline{B_t(y)}\subset B_1.
\end{equation}
Set $M=\max_{\overline{B_1}}(1-|x|)^n|h(x)|$. If $M=0$, the conclusion
is immediate. Otherwise the maximum occurs at some $y\in B_1$.
Maximality and $1-|x|\ge (1-|y|)/2$ on $B(y,(1-|y|)/2)$ imply
\begin{equation*}
 |h(x)|\le2^n|h(y)|, \forall x\in B\left(y,\frac{1-|y|}{2} \right).
\end{equation*}
Applying \eqref{eq:oscLip} with $t=(1-|y|)/2$ and setting
$C_2=2^{n+1}C_1$ yields
\begin{equation*}
 \Lip(h;B(y,(1-|y|)/4))\le\frac{C_2|h(y)|}{1-|y|}.
\end{equation*}
Choose $\tau=\min\{1/4,(2C_2)^{-1}\}$, so that $C_2\tau\le1/2$. Then
\begin{equation*}
 |h(x)|\ge |h(y)|-|h(x)-h(y)| \ge \frac{|h(y)|}{2} ,\forall x\in B(y,\tau (1-|y|)),
\end{equation*}
and integration gives
\begin{equation*}
 M\le\frac{2}{\omega_n\tau^n}\int_{B_1}|h|.
\end{equation*}
Since $1-|x|\ge1/2$ on $B_{1/2}$,
\begin{equation*}
 \sup_{B_{1/2}}|h|\le2^nM\le\frac{2^{n+1}}{\omega_n\tau^n}\int_{B_1}|h|.
\end{equation*}
A final application of \eqref{eq:oscLip} with $t=1/2$ and $y=0$ proves
\eqref{eq:L1Lip-unit}. Rescaling proves \eqref{eq:L1Lip}.
\end{proof}

\subsection{Proof of the main estimate and the \(L^q\) criterion}\label{subsec:pointwise-proof}
We now apply the preceding estimates to a sequence of comparisons on a
fixed ball. Corollary~\ref{cor:L1Lip} controls its first Lipschitz constant;
Lemma~\ref{lem:value} bounds the value error caused by one annulus;
and Proposition~\ref{prop:defect} converts this error into a Lipschitz
increment without multiplying the preceding bound. The proof first treats
smooth functions in dimensions $n>2k$. A dimensional extension and an
averaging argument then remove these restrictions.

\begin{proof}[Proof of Theorem~\ref{thm:main}]
We prove \eqref{eq:main}; the affine-invariant form will then follow by subtraction of an affine function.

\smallskip\noindent
\emph{Step 1.}
Assume first that $n>2k$ and that $u\in C^\infty\cap\Phi_k$ in a neighborhood of $\overline{B_R}$; strict admissibility is not required. By translation take $x_0=0$, and set
\[
 f=\sigma_k(D^2u),\quad \mu=f\,dx.
\]
Lemma~\ref{lem:bounded-dirichlet} provides continuous weak solutions for
all auxiliary Dirichlet problems below, whose densities are bounded and
whose boundary data are $u$ or zero. Comparison applies, and absolute
continuity ensures that spherical boundaries have zero mass.

By integration in polar coordinates, choose $s\in[R/4,R/2]$ such that
\begin{equation}\label{eq:radius}
 \int_{\partial B_s}|u|\le\frac CR\int_{B_R}|u|.
\end{equation}
By the homogeneous Dirichlet theorem \cite[Theorem~3.5]{TW}, let
$U_0\in\Phi_k(B_s)\cap C(\overline{B_s})$ be the unique $k$-harmonic
solution with boundary value $u$, and let $H$ be the ordinary harmonic
function with the same boundary value. Thus
\[
 \begin{cases}
 F_k[U_0]=0&\text{in }B_s,\\
 U_0=u&\text{on }\partial B_s,
 \end{cases}
 \qquad
 \begin{cases}
 \Delta H=0&\text{in }B_s,\\
 H=u&\text{on }\partial B_s.
 \end{cases}
\]
Comparison and the subharmonicity of $U_0$ give
\[
 u\le U_0\le H\quad\text{in }B_s.
\]
Since $|U_0|\le|u|+|H|$, and the Poisson formula gives
\[
 \int_{B_s}|H|\le\frac{s}{n}\int_{\partial B_s}|u|,
\]
\eqref{eq:radius} yields
\[
 \avint_{B_s}|U_0|\le C\avint_{B_R}|u|.
\]
By Corollary~\ref{cor:L1Lip},
\begin{equation}\label{eq:L0}
 L_0:=\Lip(U_0;B_{s/4})\le \frac CR\avint_{B_R}|u|.
\end{equation}

\medskip\noindent
\emph{Step 2.}
Let $r_j=4^{-j}s$ for $j\ge0$. On the fixed ball $B_s$, let
$U_j\in\Phi_k(B_s)\cap C(\overline{B_s})$ be the bounded-density weak
Dirichlet solution
\begin{equation}\label{eq:Uj}
 \begin{cases}
 F_k[U_j]=\mu\restr(B_s\setminus B_{r_j})&\text{in }B_s,\\
 U_j=u&\text{on }\partial B_s.
 \end{cases}
\end{equation}
The case $j=0$ agrees with the preceding definition. Since the measures
increase with $j$ and the boundary values agree, comparison gives
\[
 U_{j+1}\le U_j.
\]
Comparison with $u$ and the support of the measure in \eqref{eq:Uj} also give
\begin{equation}\label{eq:ordering}
 u\le U_{j+1}\le U_j,\quad
 F_k[U_j]=0\quad\text{in }B_{r_j}.
\end{equation}

Take a sequence $z_j\in\Phi_k(B_{2s})\cap C(\overline{B_{2s}})$ satisfying the following equation:
\begin{equation*}
 \begin{cases}
 F_k[z_j]=\mu\restr(B_{r_j}\setminus B_{r_j/4})&\text{in }B_{2s},\\
 z_j=0&\text{on }\partial B_{2s}.
 \end{cases}
\end{equation*}

Define
\[
 a_j=\left(\frac{ F_k[z_j](B_s) }{r_j^{n-k}}\right)^{1/k}.
\]

Since $z_j\le0$ by Lemma~\ref{lem:value}, 
we have
\[
 U_j+z_j\le U_j=u=U_{j+1}\text{ on }\partial B_s.
\]
Moreover, Lemma~\ref{lem:superadd}\textup{(iii)} gives
\[
 F_k[U_j+z_j]
 \ge F_k[U_j]+F_k[z_j]=F_k[U_{j+1}].
\]
The comparison principle therefore yields
\[
 U_j+z_j\le U_{j+1}\quad\text{in }B_s.
\]
Together with \eqref{eq:ordering}, this yields
\begin{equation}\label{eq:orderederror}
 0\le U_j-U_{j+1}\le-z_j.
\end{equation}
For $x\in B_{r_j/8}$, the support of $F_k[z_j]$ lies at distance
at least $r_j/8$ from $x$. Lemma~\ref{lem:value} therefore gives
\begin{equation}\label{eq:shellerror}
 -z_j(x)\le C\bigl(F_k[z_j](B_s)\bigr)^{1/k}
     \int_{r_j/8}^{8s}t^{2-n/k}\frac{dt}{t}
 \le C\bigl(F_k[z_j](B_s)\bigr)^{1/k}r_j^{2-n/k}\le C_{\rm ann}r_ja_j.
\end{equation}
The assumption $n>2k$ makes the integral controlled by its lower limit,
uniformly in $s/r_j$; fix $C_{\rm ann}=C_{\rm ann}(n,k)$ so
that the last bound holds for every $j$.

Since $U_j$ is continuous on $\overline{B_{r_j}}$ and $k$-harmonic in $B_{r_j}$,
Proposition~\ref{prop:defect} gives
\[
 \Lip(U_j;B_{r_j/2})
 \le \frac{C}{r_j}\osc_{B_{r_j}}U_j<\infty.
\]
In particular,
\[
 L_j:=\Lip(U_j;B_{r_j/4})<\infty.
\]

The function $U_{j+1}$ is $k$-harmonic in $B_{r_j/4}$. On $B_{r_j/8}$,
$U_j$ has Lipschitz constant at most $L_j$, and
\eqref{eq:orderederror}--\eqref{eq:shellerror} give
\[
 0\le U_j-U_{j+1}\le C_{\rm ann}r_ja_j.
\]
Applying \eqref{eq:LS} on $B_{r_j/8}$ with $v=U_{j+1}$, $w=U_j$ and $\delta=C_{\rm ann}r_ja_j$, we obtain
\[
 \Lip(U_{j+1};B_{r_j/16})
 \le \Lip(U_j;B_{r_j/8})
      +C_0\frac{C_{\rm ann}r_ja_j}{r_j/8}
 \le L_j+8C_0C_{\rm ann}a_j.
\]
Since $r_{j+1}/4=r_j/16$, the left-hand side is exactly $L_{j+1}$. Hence
\[
 L_{j+1}\le L_j+8C_0C_{\rm ann}a_j.
\]
Iteration gives
\begin{equation}\label{eq:iteration}
 L_m\le L_0+8C_0C_{\rm ann}\sum_{j=0}^{m-1}a_j.
\end{equation}
The coefficient one in \eqref{eq:LS} prevents geometric weights in this
sum; see Remark~\ref{rem:trudinger-comparison}.

For $r_j\le t\le2r_j$, we have $F_k[z_j](B_s)\le \mu(B_t)$ and $t^{n-k}\le2^{n-k}r_j^{n-k}$.
Therefore
\[
 \left(\frac{\mu(B_t)}{t^{n-k}}\right)^{1/k}
 \ge2^{-(n-k)/k}a_j.
\]
The intervals $[r_j,2r_j]$ are pairwise disjoint because $r_{j+1}=r_j/4$, and they lie in $(0,R]$ since $2s\le R$. Thus
\begin{equation}\label{eq:sum}
 \sum_{j=0}^\infty a_j\le C\W_k^R[\mu](0).
\end{equation}
Combining \eqref{eq:L0}, \eqref{eq:iteration}, and \eqref{eq:sum} gives the uniform bound
\begin{equation}\label{eq:Luniform}
 L_m\le \frac CR\avint_{B_R}|u|+C\W_k^R[\mu](0).
\end{equation}

\emph{Step 3.}
We pass the Lipschitz bound to rescalings converging uniformly to the
first-order affine approximation of $u$. The gradient is identified
through this affine limit.

Let $F=\|f\|_{L^\infty(B_s)}$. 
Let $z^{(j)}\in\Phi_k(B_{2s})\cap C(\overline{B_{2s}})$ be the following solution
\begin{equation*}
 \begin{cases}
 F_k[z^{(j)}]=\mu\restr B_{r_j}&\text{in }B_{2s},\\
 z^{(j)}=0&\text{on }\partial B_{2s}.
 \end{cases}
\end{equation*}
By Lemma~\ref{lem:superadd}\textup{(iii)},
\[
 F_k[U_j+z^{(j)}]
 \ge 
 F_k[u].
\]
Since $U_j+z^{(j)}\le u$ on $\partial B_s$, comparison and
\eqref{eq:ordering} yield
\[
 0\le U_j-u\le-z^{(j)}\quad\text{in }B_s.
\]
For every $x$ and $t>0$,
\[
 \mu(B_t (x)\cap B_{r_j})\le C F\min\{t^n,r_j^n\}.
\]
Using Lemma~\ref{lem:value} and splitting its integral at $r_j$, we get
\begin{equation}\label{eq:remaining}
 0\le U_j-u
 \le C F^{1/k}\left(\int_0^{r_j}t^2\frac{dt}{t}
        +r_j^{n/k}\int_{r_j}^{8s}t^{2-n/k}\frac{dt}{t}\right)\le C F^{1/k}r_j^2.
\end{equation}

Define
\[
V_j(y)=\frac{U_j(r_j y/4)-u(0)}{r_j/4}, y\in \overline{B_1}.
\]
By \eqref{eq:remaining},
\[
 \sup_{\overline{B_1}}
 \left|
V_j(y)-\frac{u(r_j y/4)-u(0)}{r_j/4}
 \right|
 \le C F^{1/k}\frac{4r_j^2}{r_j}
 \le C F^{1/k}r_j\to0.
\]
Taylor's formula gives
\[
\frac{u(r_j y/4)-u(0)}{r_j/4}
 \to Du(0)\cdot y
 \quad\text{uniformly on }\overline{B_1},
\]
and consequently
\begin{equation*}
 V_j\to V(y):=Du(0)\cdot y
 \quad\text{uniformly on }\overline{B_1}.
\end{equation*}
Also,
\[
 \Lip(V_j;B_1)
 =\Lip(U_j;B_{r_j/4})
 =L_j.
\]
Thus for every $y,z\in B_1$,
\[
 |V(y)-V(z)|
 =\lim_{j\to\infty}|V_j(y)-V_j(z)|
 \le\left(\liminf_{j\to\infty}L_j\right)|y-z|.
\]
Since $\Lip(V;B_1)=|Du(0)|$,
\[
 |Du(0)|\le\liminf_{j\to\infty}L_j.
\]
Together with \eqref{eq:Luniform}, this proves \eqref{eq:main} for
smooth functions when $n>2k$. The bound $F$ enters only the
identification of the limit, not \eqref{eq:Luniform}.

\emph{Step 4.}
If $n\le2k$, still assume that $u\in C^\infty\cap\Phi_k$ in a neighborhood
$U$ of $\overline{B_R}$. Choose $N=2k+1$ and $m=N-n$, and set
\[
 \widetilde u(x,z)=u(x),\quad (x,z)\in U\times\R^m.
\]
The Hessian has the block form
\[
 D^2_{x,z}\widetilde u=
 \begin{pmatrix}D_x^2u&0\\0&0_m\end{pmatrix}.
\]
For $1\le j\le k$,
\[
 \sigma_j^{(N)}(D^2_{x,z}\widetilde u)
 =\sigma_j^{(n)}(D_x^2u).
\]
In particular, $D_x^2u\in\overline{\Gamma_k^{\,n}}$ implies
$D^2_{x,z}\widetilde u\in\overline{\Gamma_k^{\,N}}$, and
\[
 \sigma_k^{(N)}(D^2\widetilde u)=\sigma_k^{(n)}(D^2u).
\]
Hence $\widetilde u\in C^\infty\cap\Phi_k$ on $U\times\R^m$, with $N>2k$. If $\mu=f\,dx$ and $\widetilde\mu=f(x)\,dx\,dz$, then
\begin{equation*}
 \widetilde\mu(B_t^N(x_0,0))
 =\omega_m\int_{B_t^n(x_0)}(t^2-|x-x_0|^2)^{m/2}\,d\mu(x)
 \le\omega_m t^m\mu(B_t^n(x_0)).
\end{equation*}
The factor $t^m$ cancels the dimension change in the potential, so the
$N$-dimensional potential is bounded by $C(n,k)\W_k^R[\mu](x_0)$.
Integration in the extra variables also gives
\[
 \avint_{B_R^N(x_0,0)}|\widetilde u|
 \le C(n,k)\avint_{B_R^n(x_0)}|u|,
 \quad |D\widetilde u(x_0,0)|=|Du(x_0)|.
\]
The $N$-dimensional estimate therefore proves \eqref{eq:main} for $n\le2k$.

\emph{Step 5.}
We first pass to the weak solution at a fixed averaging radius, then
shrink the averaging ball.
Let $u\in\Phi_k(\Omega)$ and fix $B_R(x_0)\Subset\Omega$. If $\W_k^R[F_k[u]](x_0)=\infty$, the estimate is immediate. We may therefore assume that this potential is finite. Choose standard local convolutions $u_m=u*\rho_{\epsilon_m}$ in a fixed neighborhood of $\overline{B_R(x_0)}$. By \eqref{eq:sobolev} and the approximation facts in Subsection~\ref{sec:preliminaries},
\[
 u_m\in C^\infty\cap\Phi_k,\quad
 Du_m\to Du\text{ in }L^1(B_R(x_0)),\quad
 \mu_m:=F_k[u_m]\to F_k[u]=\mu.
\]
Here the Hessian measures converge locally weakly.
For $y\in B_{R/4}(x_0)$, apply the smooth estimate on
$B_{R/2}(y)\Subset B_R(x_0)$ to obtain
\[
 |Du_m(y)|\le \frac CR\avint_{B_R(x_0)}|u_m|
                   +C\W_k^{R/2}[\mu_m](y).
\]
Fix $0<\rho<R/4$, average over $B_\rho(x_0)$, and let $m\to\infty$.
The $W^{1,1}$ convergence and the $L^1$ convergence of truncated
potentials from Lemma~\ref{lem:potentialcontinuity} give
\[
 \avint_{B_\rho(x_0)}|Du|
 \le \frac CR\avint_{B_R(x_0)}|u|
       +C\avint_{B_\rho(x_0)}\W_k^{R/2}[\mu](y)\,dy.
\]
Lemma~\ref{lem:averagepotential} bounds the last term by
$C\W_k^R[\mu](x_0)$. Letting $\rho\to0$ at a Lebesgue point of
$Du$ proves \eqref{eq:main}. Since subtracting an affine function
preserves $k$-convexity and the Hessian measure, application to $u-\ell$
also gives \eqref{eq:affine}.
\end{proof}

\begin{remark}
\label{rem:plaplace-iteration}
\label{rem:trudinger-comparison}
The construction in Steps~1--3 restores the measure on successive
annuli of the fixed Dirichlet problem \eqref{eq:Uj}.
The zero-order potential estimate gives
\(0\le U_j-U_{j+1}\le Cr_ja_j\) on the inner ball.
Proposition~\ref{prop:defect} converts this value error into
\[
 L_{j+1}\le L_j+Ca_j,
 \quad
 L_m\le L_0+C\sum_{j<m}a_j.
\]
Since the coefficient is one, \eqref{eq:sum} makes the increments
summable.

In the $p$-Laplace setting, the standard gradient potential argument
combines comparison with homogeneous solutions on shrinking balls and
decay of a suitable gradient excess, that is, an averaged gradient
oscillation. The recurrence is schematically
\(E_{j+1}\le\theta E_j+\text{data error}\), with \(0<\theta<1\) at
suitably chosen scales, and summability of the excess controls
successive gradient averages; see \cite{DM10,KuM14}. The
linear-potential refinement in \cite[Section~3]{KuM} uses an induction
that identifies and controls the nondegenerate regime. Here we
propagate Lipschitz bounds without decay of gradient oscillations.
For general $k$-harmonic functions with $k\ge2$ such decay is false:
\(u(x)=|x_1|\) is an example. Thus the iteration controls the size of
the gradient but gives no gradient continuity. The function average in
\eqref{eq:main} enters through the initialization \eqref{eq:L0}, not
through the iteration; affine subtraction and Poincar\'e's inequality
already give the gradient-average form \eqref{eq:gradientmean}.

A direct use of Trudinger's estimates shows the same difficulty
quantitatively. For the comparison functions and radii in Step~2, an
ordinary homogeneous Lipschitz estimate gives a recurrence with a
multiplier on \(L_j\). Assume \(n>2k\). Let \(C_T\) denote a constant
in the centered homogeneous estimate
\[
 |Dh(y)|\le\frac{C_T}{\rho}\osc_{B_\rho(y)}h
\]
of \cite[Theorem~3.1]{Tru97} for smooth homogeneous solutions, and use
its Lipschitz consequence for continuous weak solutions by
approximation. For \(y\in B_{r_j/16}\), the ball
\(B_{r_j/16}(y)\) lies in \(B_{r_j/8}\), so
\eqref{eq:orderederror}--\eqref{eq:shellerror} give
\[
 \osc_{B_{r_j/16}(y)}U_{j+1}
 \le \frac{r_j}{8}L_j+Kr_ja_j
\]
with a fixed \(K=K(n,k)\). Applying the centered estimate directly
gives only
\[
 L_{j+1}\le 2C_TL_j+16C_TK a_j.
\]
For \(k\ge2\), smooth convex approximations of \((x_1+a)_+\) on \(B_1\),
with \(0<a<1\), have Hessians of rank at most one and are
\(k\)-harmonic. Applying the centered estimate to these functions gives
\(C_T\ge(1+a)^{-1}\), hence \(C_T\ge1\) as \(a\to0\). Iterating the
preceding bound therefore produces geometric weights:
\[
 L_m\le (2C_T)^mL_0
       +16C_TK\sum_{j=0}^{m-1}(2C_T)^{m-1-j}a_j.
\]
The summability in \eqref{eq:sum} does not give a bound uniform in
\(m\), because \(L_j\) is multiplied repeatedly. Trudinger's
H\"older estimate \cite[Theorem~4.1]{Tru97} has \(\alpha<1\). Its proof
uses the small factor \((\sigma/\tau)^{1-\alpha}\) to absorb the
homogeneous oscillation term; see \cite[(4.5)--(4.9)]{Tru97}. This
absorption is unavailable at \(\alpha=1\). The comparison \eqref{eq:LS}
gives instead the coefficient-one estimate leading to the additive
bound \eqref{eq:iteration}.
\end{remark}

The first consequence of the pointwise estimate follows by bounding the
potential through H\"older's inequality. Its scale exponent is positive
precisely when $q>n/k$.

\begin{proof}[Proof of Corollary~\ref{cor:lqlipschitz}]
Write $N=\|f\|_{L^q(B_R(x_0))}$ and $\omega_n=|B_1|$. For $x\in B_{R/2}(x_0)$ and $0<t\le R/2$, H\"older's inequality gives
\[
 F_k[u](B_t(x))
 =\int_{B_t(x)}f\,dy
 \le\omega_n^{1-1/q}N t^{n(1-1/q)},
\]
with the same interpretation when $q=\infty$. Consequently,
\begin{equation}\label{eq:lqpotentialbound}
 \sup_{x\in B_{R/2}(x_0)}\W_k^{R/2}[F_k[u]](x)
 \le C(n,k)N^{1/k}\int_0^{R/2}t^{-n/(kq)}\,dt
 \le\frac{C(n,k)}{\alpha_q}R^{\alpha_q}N^{1/k}.
\end{equation}
Apply \eqref{eq:affine} on $B_{R/2}(x)\subset B_R(x_0)$ at Lebesgue
points $x\in B_{R/2}(x_0)$ of $Du$; the average term is bounded by $2^n$ times the corresponding
average on $B_R(x_0)$. Taking the essential supremum and using
\eqref{eq:lqpotentialbound} proves \eqref{eq:lqlipschitz}. Since
$u\in W^{1,1}_{\mathrm{loc}}(\Omega)$, the locally bounded weak gradient
gives \eqref{eq:lqcriterion} and a locally Lipschitz representative.
\end{proof}

\section{Equations with gradient growth}\label{sec:gradient-growth}\label{sharp:gradient-growth}
This section proves Theorem~\ref{growth:lipschitz} and establishes
the sharpness of its exponent in the full weak class. The basic
Sobolev integrability of $k$-convex functions starts an iteration,
which closes once the density exponent exceeds $n/k$.

\begin{proof}[Proof of Theorem~\ref{growth:lipschitz}]
All estimates are applied on nested relatively compact balls. We first
record the density mapping needed below:
\begin{equation}\label{growth:density-mapping}
 f\in L^t_{\mathrm{loc}}(\Omega),\quad 1<t<n/k
 \quad\Longrightarrow\quad
 Du\in L^{nkt/(n-kt)}_{\mathrm{loc}}(\Omega).
\end{equation}
Indeed, restrict $f$ to a fixed ball and apply \eqref{eq:main} on
balls of a fixed smaller radius centered at its interior points.
The mean term is uniformly bounded, while the nonlinear potential
mapping in \cite[Proposition~2.8 and Theorem~3.1(i)]{Cianchi}, with
parameters $k/(k+1)$ and $k+1$, maps $L^t$ into
$L^{nkt/(n-kt)}$. This gives \eqref{growth:density-mapping}.
By \eqref{eq:sobolev}, $Du\in L^r_{\mathrm{loc}}$ for
$1\le r<p_*$. Since $p<p_*$ and $n(p_*-k)/k=p_*$, choose
\begin{equation}\label{growth:initial-exponent}
 \max\left\{1,p,\frac{n(p-k)}k\right\}<r_0<p_*.
\end{equation}

Our aim is to reach $r>np/k$, so that $|Du|^p$ has density
exponent $r/p>n/k$. If $r_0$ already exceeds this threshold,
no iteration is needed. Otherwise we use the following improvement.
Suppose $Du\in L^r_{\mathrm{loc}}$ with $r\ge r_0$ and $r<np/k$.
Then $1<r/p<n/k<q$, so \eqref{growth:equation} gives
$f\in L^{r/p}_{\mathrm{loc}}$. Applying
\eqref{growth:density-mapping} with $t=r/p$ yields
\begin{equation*}
 Du\in L^{r_+}_{\mathrm{loc}},
 \quad
 r_+=\frac{nkr}{np-kr},\quad
 \frac1{r_+}=\frac pk\frac1r-\frac1n.
\end{equation*}
Moreover,
\begin{equation*}
 r_+>r\quad\Longleftrightarrow\quad r>\frac{n(p-k)}k,
\end{equation*}
which holds by \eqref{growth:initial-exponent}.

The iteration reaches $r\ge np/k$ in finitely many steps. Otherwise
$s_j=1/r_j>0$ would satisfy $s_{j+1}=(p/k)s_j-1/n$ indefinitely.
If $p\le k$, each step decreases $s_j$ by at least $1/n$.
If $p>k$, let $s_*=k/[n(p-k)]$; then
$s_j-s_*=(p/k)^j(s_0-s_*)$, with $s_0<s_*$ by
\eqref{growth:initial-exponent}. Both cases eventually give a negative
$s_j$, a contradiction. Only finitely many improvements are
therefore needed; the intermediate balls may all be chosen to
contain a prescribed compact subset of $\Omega$.

If $r=np/k$, we have $f\in L^{n/k}_{\mathrm{loc}}$.
For any $1<t<n/k$, \eqref{growth:density-mapping} then gives
$Du\in L^{nkt/(n-kt)}_{\mathrm{loc}}$. Choosing $t$ sufficiently
close to $n/k$ makes this exponent exceed $np/k$; relabel the
larger exponent as $r$.
Thus in all cases
\[
 f\in L^{\min\{r/p,q\}}_{\mathrm{loc}}(\Omega),\quad
 \min\{r/p,q\}>n/k.
\]
Corollary~\ref{cor:lqlipschitz} completes the proof.
\end{proof}

The conclusion concerns gradient size, not gradient continuity.
For instance, $u(x)=|x_1|$ satisfies the growth assumption with
$f=h=0$ for $k\ge2$ but is not $C^1$; see
Remark~\ref{rem:laplacian-density-example}.
The strict inequality in \eqref{growth:range} is essential in the
full weak class: Proposition~\ref{growth:critical-example} constructs
a solution of $F_k[u]=|Du|^{p_*}$ with $Du\in L^{p_*}$ and an
unbounded gradient, as the following construction shows.
The density must now equal the stated power of the gradient. The
following radial example verifies this equality for the full Hessian measure.

\begin{proposition}\label{growth:critical-example}
For every $2\le k<n$, let $p_*=nk/(n-k)$. There is a radial function
\begin{equation*}
 u\in\Phi_k(B_1)\cap W^{1,p_*}(B_1),\quad
 F_k[u]=|Du|^{p_*}\quad\text{in }B_1,
\end{equation*}
which is smooth and strictly $k$-convex on $B_1\setminus\{0\}$ and
satisfies $|Du(x)|\to\infty$ as $x\to0$, $x\ne0$.
Its full Hessian measure has no singular part.
If $n<2k$, the function is also continuous at zero.
\end{proposition}

\begin{proof}
Set
\begin{equation*}
 a=\frac{n-k}{k},\quad b=\frac{n-k}{k^2},\quad
 c=\left[\binom{n-1}{k-1}b\right]^{1/(p_*-k)},\quad
 L(r)=\log\frac er,
\end{equation*}
and define
\begin{equation*}
 w(r)=c r^{-a}L(r)^{-b},\quad
 U(r)=-\int_r^1 w(t)\,dt,\quad u(x)=U(|x|).
\end{equation*}
Give $u(0)$ its radial limit, finite if $n<2k$ and equal to $-\infty$
otherwise. Since $U'=w>0$ and $rw'/w=-a+b/L$, the radial eigenvalue
formula gives, for $1\le j\le k$,
\begin{equation}\label{growth:radial-admissibility}
 \sigma_j(D^2u)=\binom{n-1}{j-1}
       \left(\frac{w(r)}r\right)^j
       \left(\frac nj-\frac nk+\frac b{L(r)}\right)>0.
\end{equation}
The identities $ap_*=n$ and $bp_*=n/k=kb+1$, together with the
choice of $c$, yield
\begin{equation*}
 \sigma_k(D^2u)
 =\frac1k\binom{n-1}{k-1}r^{1-n}\big(r^{n-k}w(r)^k\big)'
   =w(r)^{p_*}=|Du|^{p_*},\quad r>0.
\end{equation*}

Although $w(r)\to\infty$,
\begin{equation}\label{growth:critical-integrability}
 \int_{B_1}|Du|^{p_*}
 =n\omega_n c^{p_*}\int_0^1\frac{dr}{rL(r)^{n/k}}<\infty.
\end{equation}
Since $L(t)^{-b}\le1$, the defining integral gives
\[
 |U(r)|\le C
 \begin{cases}
  1,&a<1,\\
  \log(e/r),&a=1,\\
  r^{1-a},&a>1.
 \end{cases}
\]
In the last case $n+(1-a)p_*=p_*>0$, so these bounds give
$u\in L^{p_*}(B_1)$ in every case. They also imply
$r^{n-1}|U(r)|\to0$. Integration by parts on punctured balls,
together with \eqref{growth:critical-integrability}, therefore
identifies $Du=w(r)x/r$ as the weak gradient and proves
$u\in W^{1,p_*}(B_1)$.

For the full measure, $0<\varepsilon<1$, keep $U_\varepsilon=U$ on $[\varepsilon,1)$ and set
\[
 U_\varepsilon(r)=U(\varepsilon)
       +\frac{w(\varepsilon)}{2\varepsilon}(r^2-\varepsilon^2),
 \quad 0\le r\le\varepsilon.
\]
The matching first derivatives give a locally $C^{1,1}$ radial
function $u_\varepsilon$ whose Hessian is in $\Gamma_k$ almost
everywhere: it is a positive multiple of $I$ inside $B_\varepsilon$,
and \eqref{growth:radial-admissibility} applies outside.
For fixed $\varepsilon$, convolution preserves admissibility and the
mollified Hessians converge in $L^k_{\mathrm{loc}}$. Their $k$th
symmetric functions therefore converge in $L^1_{\mathrm{loc}}$;
\eqref{eq:weakcont} identifies
$F_k[u_\varepsilon]=\sigma_k(D^2u_\varepsilon)\,dx$, with no measure
on the joining sphere. Its inner mass is
\begin{equation}\label{growth:zero-mass}
 F_k[u_\varepsilon](B_\varepsilon)
 =\omega_n\binom nk\varepsilon^{n-k}w(\varepsilon)^k
 =\omega_n\binom nk c^k L(\varepsilon)^{-kb}
   \to0\text{ as }\varepsilon\to0.
\end{equation}

Outside this ball the density is $|Du|^{p_*}$. Moreover,
\[
 \int_{B_\varepsilon}|u_\varepsilon-u|
 \le C\varepsilon^n\bigl(|U(\varepsilon)|
                   +\varepsilon w(\varepsilon)\bigr)
       +\int_{B_\varepsilon}|u|\to0
\]
by the preceding bounds. Since $u_\varepsilon=u$ off
$B_\varepsilon$, this proves local $L^1$ convergence.
Closure of $\Phi_k$ and weak
continuity \eqref{eq:weakcont}, together with
\eqref{growth:critical-integrability}--\eqref{growth:zero-mass}, show
that $F_k[u]=|Du|^{p_*}\,dx$ with no mass at zero.
Finally, $n<2k$ is exactly $a<1$, which makes the defining integral
finite at zero and gives the stated continuity.
\end{proof}

\section{Regularity consequences}\label{sec:potential-scale}\label{sec:data-spaces}\label{sec:consequences}

Throughout this section, we assume $2\le k<n$ and
$u\in\Phi_k(\Omega)$, and write $\mu=F_k[u]$.
We apply Theorem~\ref{thm:main} to obtain integrability and
boundedness estimates for $Du$ and oscillation estimates for $u$.
Sharpness examples are given in the
\hyperref[sec:examples]{appendix}.

\subsection{Gradient integrability and boundedness}\label{data:critical-subsection}\label{sharp:zero-measure}
Fix $B=B_R(z)\Subset\Omega$ and $a\in\R^n$, and set
$E=B_{R/2}(z)$ and $A=\avint_B|Du-a|$. Let $\lambda$ denote
the extension of $\mu\restr B$ by zero to $\R^n$.
For almost every $x\in E$, we apply \eqref{eq:gradientmean} on
$B_{R/2}(x)\subset B$. The average is bounded by $2^n A$ and
the truncated potential by $\W_k^\infty[\lambda](x)$, giving
\begin{equation}\label{data:local-reduction}
 |Du(x)-a|\le CA+C\W_k^\infty[\lambda](x)
 \quad\text{for almost every }x\in E.
\end{equation}
Using \eqref{eq:affine} instead, we may replace $A$ throughout
by $R^{-1}\avint_B|u-\ell|$ for any affine function $\ell$
with $D\ell=a$.

For $0<q<\infty$ and $0<s<\infty$, we use the Lorentz quasi-norms
\begin{equation*}
 \|F\|_{L^{q,s}(G)}
 =\left(\int_0^{|G|}[t^{1/q}F^*(t)]^s\frac{dt}{t}\right)^{1/s},
 \quad
 \|F\|_{L^{q,\infty}(G)}
 =\sup_{0<t<|G|}t^{1/q}F^*(t),
\end{equation*}
where $F^*$ is the decreasing rearrangement of $|F|$ extended
by zero outside $G$. The upper limits are infinite when
$|G|=\infty$. These definitions apply also for $s<1$, and
$L^{q,q}(G)=L^q(G)$.

We use the mapping estimates in
\cite[Proposition~2.8 and Theorem~3.1]{Cianchi} with
$\gamma=k/(k+1)$ and $p=k+1$. These parameters satisfy
$\gamma p=k<n$ and $p-1=k$. In particular,
\cite[(4.4)]{Cianchi} gives, for $F\ge0$,
\begin{equation}\label{data:iterated-riesz-estimate}
 \W_k^\infty[F\,dx]\le C I_\gamma[(I_\gamma F)^{1/k}],
 \quad
 I_\gamma F(x)=\int_{\R^n}|x-y|^{\gamma-n}F(y)\,dy.
\end{equation}

\begin{proposition}\label{data:integrability}
Let $p_*=nk/(n-k)$. With the notation above, the following estimates hold.
\begin{enumerate}
\renewcommand{\labelenumi}{\textup{(\roman{enumi})}}
\renewcommand{\theenumi}{\thetheorem(\roman{enumi})}
\item For arbitrary Hessian measure,
\begin{equation}\label{data:measure-weak-estimate}
 \|Du-a\|_{L^{p_*,\infty}(E)}
 \le CR^{n/p_*}A+CF_k[u](B)^{1/k},
\end{equation}
where $C=C(n,k)$.

\item If $F_k[u]=f$ on $B$, $1<q<n/k$, $0<s\le\infty$,
and $p=nkq/(n-kq)$, then
\begin{equation}\label{data:lorentz-estimate}
 \|Du-a\|_{L^{p,ks}(E)}
 \le CR^{n/p}A+C\|f\|_{L^{q,s}(B)}^{1/k},
\end{equation}
where $k\infty=\infty$ and $C=C(n,k,q,s)$.

\item\label{data:critical-lipschitz}
If $F_k[u]=f$ on $B$ and $f\in L^{n/k,1/k}(B)$, then
\begin{equation}\label{data:critical-lipschitz-estimate}
 \|Du-a\|_{L^\infty(E)}
 \le CA+C\|f\|_{L^{n/k,1/k}(B)}^{1/k},
\end{equation}
where $C=C(n,k)$. If this density condition holds locally throughout
$\Omega$, then $u$ is locally Lipschitz and, for every compact
$K\Subset\Omega$,
\begin{equation}\label{data:critical-uniform-tails}
 \lim_{\rho\to0}\sup_{x\in K}\W_k^\rho[f\,dx](x)=0.
\end{equation}

\item If $F_k[u]=f$ on $B$ and
$f\in L^{q,\infty}(B)$ for $n/k<q\le\infty$, let
$\eta_q=1-n/(kq)>0$, with $\eta_\infty=1$ and
$L^{\infty,\infty}=L^\infty$. Then
\[
 \|Du-a\|_{L^\infty(E)}
 \le CA+\frac{C}{\eta_q}R^{\eta_q}
       \|f\|_{L^{q,\infty}(B)}^{1/k}.
\]
Moreover, for every $\rho>0$,
\begin{equation}\label{data:supercritical-tail}
 \sup_{x\in\R^n}\W_k^\rho[f\mathbf1_B\,dx](x)
 \le\frac{C}{\eta_q}\rho^{\eta_q}
       \|f\|_{L^{q,\infty}(B)}^{1/k},
\end{equation}
where $C=C(n,k)$. In particular, the small-scale potentials vanish
uniformly on compact subsets of $B$.
\end{enumerate}
\end{proposition}

\begin{proof}
\emph{(i)} Let $m=\lambda(\R^n)$ and let $M\lambda$ be the centered
maximal function of the measure. Splitting the Wolff integral at
$\delta>0$ gives
\[
 \W_k^\infty[\lambda](x)
 \le C(M\lambda(x))^{1/k}\delta
       +Cm^{1/k}\delta^{1-n/k}.
\]
For $0<M\lambda(x)<\infty$, take $\delta=(m/M\lambda(x))^{1/n}$.
Then
$\W_k^\infty[\lambda]\le Cm^{1/n}(M\lambda)^{1/k-1/n}$.
Since $(1/k-1/n)^{-1}=p_*$, the weak $(1,1)$ maximal estimate yields
\[
 |\{x:\W_k^\infty[\lambda](x)>\tau\}|
 \le C\left(\frac{m^{1/k}}{\tau}\right)^{p_*},\quad \tau>0.
\]
Thus $\|\W_k^\infty[\lambda]\|_{L^{p_*,\infty}}\le Cm^{1/k}$;
the case $m=0$ is immediate. Combining this with
\eqref{data:local-reduction} proves \eqref{data:measure-weak-estimate}.

\emph{(ii)} Apply \cite[Theorem~3.1(i)]{Cianchi} to
\eqref{data:iterated-riesz-estimate}, with $F=f\mathbf1_B$.
The target space is $L^{nkq/(n-kq),ks}$ for all $0<s\le\infty$,
including $ks<1$. Equation~\eqref{data:local-reduction} and
$\|A\|_{L^{p,ks}(E)}=CR^{n/p}A$ give
\eqref{data:lorentz-estimate}.

\emph{(iii)} By \cite[Theorem~3.1(iv)]{Cianchi} and
\eqref{data:iterated-riesz-estimate},
\begin{equation*}
 \sup_{x\in\R^n}\W_k^\infty[F\,dx](x)
 \le C\|F\|_{L^{n/k,1/k}(\R^n)}^{1/k}.
\end{equation*}
The essential bound holds at every point because the potential is
lower semicontinuous. Taking $F=f\mathbf1_B$ and using
\eqref{data:local-reduction} proves
\eqref{data:critical-lipschitz-estimate}.

For the uniform tails, choose $K\Subset U\Subset\Omega$ and extend
$F=f\mathbf1_U$ by zero. For small $\rho$ and $x\in K$,
\begin{equation*}
 \W_k^\rho[f\,dx](x)
 \le\W_k^\infty[F\mathbf1_{B_\rho(x)}\,dx](x)
 \le C\|F\mathbf1_{B_\rho(x)}\|_{L^{n/k,1/k}}^{1/k}.
\end{equation*}
For a measurable set $G$,
$(F\mathbf1_G)^*(t)\le F^*(t)\mathbf1_{(0,|G|)}(t)$, and hence
\[
 \|F\mathbf1_G\|_{L^{n/k,1/k}}^{1/k}
 \le\int_0^{|G|}[t^{k/n}F^*(t)]^{1/k}\frac{dt}{t}
 \to0\quad(|G|\to0).
\]
The bound depends only on $|G|$. Taking $G=B_\rho(x)$ proves
\eqref{data:critical-uniform-tails}.
The bounded weak gradient gives a locally Lipschitz representative,
which agrees with the canonical subharmonic representative of $u$.

\emph{(iv)} For finite $q$, the weak-$L^q$ estimate gives
\[
 \int_{B_t(x)}f\mathbf1_B\,dx
 \le\frac{q}{q-1}|B_t|^{1-1/q}
       \|f\|_{L^{q,\infty}(B)}.
\]
Since $q>n/k>1$, its coefficient is bounded in terms of $n,k$.
Substitution in the Wolff integral gives
\eqref{data:supercritical-tail}, because
$\int_0^\rho t^{\eta_q}\,dt/t=\rho^{\eta_q}/\eta_q$.
For $q=\infty$, use
$\int_{B_t(x)}f\mathbf1_B\,dx\le |B_t|\|f\|_\infty$ instead.
Applying \eqref{eq:gradientmean} on
$B_{R/2}(x)\subset B$ for almost every $x\in E$ proves the gradient
bound; its average term is at most $2^nA$.
\end{proof}

\begin{remark}
\label{rem:laplacian-density-example}\label{obs:positive-exponents}
The space \(L^{p_*,\infty}\) in \eqref{data:measure-weak-estimate}
cannot in general be replaced by \(L^{p_*}\), where
\(p_*=nk/(n-k)\). The radial \(k\)-Hessian fundamental solution is a
suitably signed multiple of \(|x|^{2-n/k}\) when \(n\ne2k\), and a
positive multiple of \(\log|x|\) when \(n=2k\); see \cite{TW}. Its
Hessian measure is \(c\delta_0\) with \(c>0\), and
\(|Du(x)|\asymp |x|^{1-n/k}\), so \(Du\in L^{p_*,\infty}\) but
\(Du\notin L^{p_*}\) near the origin. The exponent \(1/k\) in
Proposition~\ref{data:critical-lipschitz} is sharp as well: the example
in Proposition~\ref{prop:criticaldensity} has unbounded gradient, while
its Hessian density satisfies
\[
f(x)\asymp |x|^{-k}\bigl[\log(4/|x|)\bigr]^{-k},
\quad
f\in L^{n/k,s}(B_1)\quad\text{for every }s>1/k.
\]

The potential can also vanish identically while the gradient is
discontinuous. For \(u(x)=|x_1|\), the smooth convex approximations
\(u_\varepsilon(x)=\sqrt{x_1^2+\varepsilon^2}\) have rank-one
Hessians, so weak continuity gives \(F_k[u]=0\) for \(k\ge2\), whereas
\(Du=\operatorname{sgn}(x_1)e_1\) almost everywhere. On every ball
\(B_\rho\) centered at the origin,
\[
\avint_{B_\rho}|Du-(Du)_{B_\rho}|=1,
\]
so \(Du\) is not in VMO. Gradient BMO can fail even under vanishing
critical Morrey density; see Proposition~\ref{obs:bmo-limits}.


\end{remark}

\subsection{Averaged estimates and moduli of continuity}\label{sharp:logarithmic}
The averaged estimates below hold at every center.
For $B_R(x)\subset\Omega$, $\beta\ge0$, and $0\le\alpha\le1$, set
\begin{equation*}
 M_{\beta,R}F(x)
 =\sup_{0<r\le R}r^\beta\avint_{B_r(x)}|F|,
 \quad
 M^\#_{\alpha,R}v(x)
 =\sup_{0<r\le R}r^{-\alpha}
      \avint_{B_r(x)}|v-(v)_{B_r(x)}|.
\end{equation*}
For $0\le\alpha\le1$ and $\sigma\ge0$, we also define
\begin{equation*}
 P_\alpha^R[\mu](x)
 =\int_0^R
   \left(\frac{\mu(B_t(x))}{t^{n-2k+k\alpha}}\right)^{1/k}
   \frac{dt}{t},
 \quad
 Q_{\sigma,R}\mu(x)
 =\sup_{0<t\le R}\frac{\mu(B_t(x))}{t^{n-\sigma}}.
\end{equation*}
With this notation, $P_1^R=\W_k^R$. We use $Q_{\sigma,R}\mu(x)$
through the bound
$\mu(B_t(x))\le Q_{\sigma,R}\mu(x)t^{n-\sigma}$ for $0<t\le R$,
which holds also for $\sigma\ge n$.

\begin{proposition}\label{pot:uniform}
Let $B_R(x)\Subset\Omega$, let $\ell$ be affine with $D\ell=a$, and set
$A=\avint_{B_R(x)}|Du-a|$. Then
\begin{align}
 M^\#_{\alpha,R}(u-\ell)(x)+M_{1-\alpha,R}(Du-a)(x)
 &\le CP_\alpha^R[\mu](x)+CR^{1-\alpha}A,
 \quad 0\le\alpha\le1,\notag\\
 M^\#_{\alpha,R}(u-\ell)(x)+M_{1-\alpha,R}(Du-a)(x)
 &\le\frac{C}{1-\alpha}
       [Q_{k(2-\alpha),R}\mu(x)]^{1/k}+CR^{1-\alpha}A,
 \quad 0\le\alpha<1,\label{pot:maximal-estimate}
\end{align}
where $C=C(n,k)$ is independent of $\alpha$. Infinite right-hand sides
are allowed.
\end{proposition}
\begin{proof}
The sharp maximal term follows from the averaged gradient term by
Poincar\'e's inequality. Indeed, since $u\in W^{1,1}_{\mathrm{loc}}$,
for every $0<r\le R$ we have
\[
 r^{-\alpha}\avint_{B_r(x)}
   |(u-\ell)-(u-\ell)_{B_r(x)}|
 \le C r^{1-\alpha}\avint_{B_r(x)}|Du-a|.
\]
Taking the supremum in $r$ gives
\[
 M^\#_{\alpha,R}(u-\ell)(x)
 \le C M_{1-\alpha,R}(Du-a)(x).
\]
It therefore suffices to establish the two stated estimates for
$M_{1-\alpha,R}(Du-a)(x)$.

Fix $0<r\le R/4$. For almost every $y\in B_r(x)$, apply
\eqref{eq:gradientmean} on $B_{R/2}(y)\subset B_R(x)$.
Its average term is at most $2^n A$. We average the potential term
over $B_r(x)$ and split its radii at $r$. For $0<t<r$, Fubini's
theorem gives
\[
 \avint_{B_r(x)}\mu(B_t(y))\,dy
 =\frac1{|B_r|}\int_{B_{r+t}(x)}
    |B_r(x)\cap B_t(z)|\,d\mu(z)
 \le C(t/r)^n\mu(B_{2r}(x)).
\]
By concavity of the $k$th root,
\[
 \avint_{B_r(x)}\mu(B_t(y))^{1/k}\,dy
 \le C(t/r)^{n/k}\mu(B_{2r}(x))^{1/k}.
\]
Integrating against $t^{1-n/k}\,dt/t$ over $(0,r)$ gives
$C\mu(B_{2r}(x))^{1/k}r^{1-n/k}$. For $r\le t\le R/2$,
use $B_t(y)\subset B_{2t}(x)$ and then change variables from $t$
to $2t$. Thus
\begin{equation}\label{pot:average-splitting}
 \avint_{B_r(x)}|Du-a|
 \le CA+C\mu(B_{2r}(x))^{1/k}r^{1-n/k}
       +C\int_{2r}^R\mu(B_t(x))^{1/k}t^{1-n/k}\frac{dt}{t}.
\end{equation}
All integrands are nonnegative, so these steps remain valid when a
potential is infinite.

Multiply \eqref{pot:average-splitting} by $r^{1-\alpha}$.
For $0\le\alpha\le1$, monotonicity of the ball mass gives
\[
 \mu(B_{2r}(x))^{1/k}r^{2-\alpha-n/k}
 \le C\int_{2r}^{4r}
     \left(\frac{\mu(B_t(x))}{t^{n-2k+k\alpha}}\right)^{1/k}
     \frac{dt}{t}.
\]
The constant is uniform in $\alpha$, because $t/r\in[2,4]$.
For the final term in \eqref{pot:average-splitting}, use
$r^{1-\alpha}\le t^{1-\alpha}$. Both measure terms are consequently
bounded by $CP_\alpha^R[\mu](x)$, and the average term is at most
$CR^{1-\alpha}A$.

For the second estimate, let $0\le\alpha<1$ and set
$Q=Q_{k(2-\alpha),R}\mu(x)$. We may assume $Q<\infty$.
Then $\mu(B_t(x))\le Qt^{n-2k+k\alpha}$ for $t\le R$.
Substitution into \eqref{pot:average-splitting} yields
\[
 r^{1-\alpha}\avint_{B_r(x)}|Du-a|
 \le CR^{1-\alpha}A+
 CQ^{1/k}\left(1+r^{1-\alpha}
                      \int_{2r}^R t^{\alpha-1}\frac{dt}{t}\right)
 \le CR^{1-\alpha}A+\frac{C}{1-\alpha}Q^{1/k}.
\]
Finally, for $R/4<r\le R$,
\[
 r^{1-\alpha}\avint_{B_r(x)}|Du-a|
 \le 4^nR^{1-\alpha}A.
\]
Taking the supremum over all radii proves both gradient maximal
estimates. The Poincar\'e reduction at the beginning then supplies
the sharp maximal term in each inequality.
\end{proof}

\begin{corollary}\label{pot:moduli}
Let $B_{2R}(z)\Subset\Omega$, let $\ell$ be affine, and set
$A=\avint_{B_{2R}(z)}|Du-D\ell|$.
\begin{enumerate}
\renewcommand{\labelenumi}{\textup{(\roman{enumi})}}
\renewcommand{\theenumi}{\thetheorem(\roman{enumi})}
\item Suppose that, for some $0\le\alpha<1$,
\begin{equation}\label{pot:morrey-assumption}
 \mu(B_t(x))\le Mt^{n-2k+k\alpha}
 \quad(x\in B_R(z),\ 0<t\le R).
\end{equation}
Then
\begin{align}
 [u-\ell]_{C^{0,\alpha}(B_{R/4}(z))}
 &\le\frac{C}{\alpha}
       \left(\frac{M^{1/k}}{1-\alpha}+R^{1-\alpha}A\right),
 \quad 0<\alpha<1,\label{pot:holder-estimate}\\
 [u-\ell]_{\mathrm{BMO}(B_{R/4}(z))}
 &\le C(M^{1/k}+RA),\quad\alpha=0.\label{pot:bmo-function}
\end{align}
Here BMO uses mean oscillations on subballs.

\item\label{pot:log-lipschitz}
Suppose instead that
\begin{equation}\label{pot:critical-growth}
 \mu(B_t(x))\le Mt^{n-k}
 \quad(x\in B_R(z),\ 0<t\le R).
\end{equation}
Then, for distinct $x,y\in B_{R/4}(z)$,
\begin{equation}\label{pot:log-modulus}
 |(u-\ell)(x)-(u-\ell)(y)|
 \le C|x-y|\left[A+M^{1/k}\log\frac{eR}{|x-y|}\right].
\end{equation}
\end{enumerate}
In both parts $C=C(n,k)$.
\end{corollary}

\begin{proof}
\emph{(i)} For every $x\in B_R(z)$, one has
$B_R(x)\subset B_{2R}(z)$ and
$\avint_{B_R(x)}|Du-D\ell|\le2^nA$. By
\eqref{pot:maximal-estimate} and \eqref{pot:morrey-assumption},
\[
 M^\#_{\alpha,R}(u-\ell)(x)
 \le\frac{CM^{1/k}}{1-\alpha}+CR^{1-\alpha}A.
\]
For $0<\alpha<1$, summing differences of averages on
$B_{2^{-j}r}(x)$ gives the factor $(1-2^{-\alpha})^{-1}$.
Comparing averages at two centers on the scale of their distance then
proves \eqref{pot:holder-estimate}, since
$(1-2^{-\alpha})^{-1}\le C/\alpha$. The resulting continuous
representative agrees with the canonical subharmonic representative.
For $\alpha=0$, the same maximal estimate directly bounds mean
oscillations on every subball of $B_{R/4}(z)$ and proves
\eqref{pot:bmo-function}.

\emph{(ii)} For $0<\beta\le1/2$, \eqref{pot:critical-growth}
implies $Q_{k(1+\beta),R}\mu(x)\le MR^{k\beta}$ uniformly for
$x\in B_R(z)$. Apply \eqref{pot:maximal-estimate} with
$\alpha=1-\beta$ and the Campanato argument in part~\textup{(i)}.
Its constant is bounded because $\alpha\ge1/2$. With $d=|x-y|$,
\[
 |(u-\ell)(x)-(u-\ell)(y)|
 \le Cd^{1-\beta}R^\beta(A+M^{1/k}/\beta).
\]
Choose $\beta=\min\{1/2,1/\log(eR/d)\}$. Then
$d^{1-\beta}R^\beta=d(R/d)^\beta\le ed$ and
$1/\beta\le2\log(eR/d)$, which give \eqref{pot:log-modulus}.
\end{proof}

Suppose that $\mu=f\,dx$ on $B_{2R}(z)$. For $1<q<\infty$
and every measurable set $G\subset B_{2R}(z)$, we have
\[
 \int_G f\,dx
 \le \frac{q}{q-1}|G|^{1-1/q}
       \|f\|_{L^{q,\infty}(B_{2R}(z))}.
\]
Taking $G=B_t(x)$ gives
$\mu(B_t(x))\le C(n,q)\|f\|_{L^{q,\infty}(B_{2R}(z))}t^{n-n/q}$
for $x\in B_R(z)$ and $0<t\le R$.
If $\max\{1,n/(2k)\}<q<n/k$, then
$\alpha=2-n/(kq)\in(0,1)$ and $n-n/q=n-2k+k\alpha$.
Corollary~\ref{pot:moduli} therefore yields
\begin{equation}\label{data:weak-holder}
 [u-\ell]_{C^{0,\alpha}(B_{R/4}(z))}
 \le C(n,k,q)\left[
 \|f\|_{L^{q,\infty}(B_{2R}(z))}^{1/k}
 +R^{1-\alpha}A\right].
\end{equation}
At $q=n/k$, the same calculation gives the growth exponent
$n-k$. Applying \eqref{pot:log-modulus}, we obtain
\begin{equation*}
 |(u-\ell)(x)-(u-\ell)(y)|
 \le C|x-y|\left[
 A+\|f\|_{L^{n/k,\infty}(B_{2R}(z))}^{1/k}
       \log\frac{eR}{|x-y|}\right]
\end{equation*}
for distinct $x,y\in B_{R/4}(z)$, where $C=C(n,k)$.

\begin{remark}
The H\"older exponent in \eqref{data:weak-holder} is optimal.
Let $\max\{1,n/(2k)\}<q<n/k$,
$\alpha=2-n/(kq)$, and $u(x)=|x|^\alpha$. For $x\ne0$,
\[
 \sigma_k(D^2u)
 =\frac1k\binom{n-1}{k-1}
   \alpha^k(n-2k+k\alpha)|x|^{k\alpha-2k}
 \in L^{q,\infty}(B_1).
\]
The radial flux tends to zero at the origin, so the full
Hessian measure is given by this density. Since
$u(re)-u(0)=r^\alpha$ for every unit vector $e$, the H\"older
exponent at the origin is exactly $\alpha$.

The logarithmic factor in \eqref{pot:log-modulus} is also sharp,
as shown by the example
\[
 u(x)=x_1\log\frac4{|x|}
       +\frac{6n}{n-k}|x|
\]
in Proposition~\ref{obs:loglip-sharp}.

The folding construction in Proposition~\ref{obs:bmo-limits}
gives a function with vanishing critical Morrey density and
$Du\notin\mathrm{BMO}_{\mathrm{loc}}$.
For $k\ge3$, smoothing the fold at a variable scale gives
a density in $L^{n/k}$ while preserving the failure of
gradient BMO.
\end{remark}

\renewcommand{\theHsubsection}{appendix.\arabic{subsection}}
\renewcommand{\theHtheorem}{appendix.\arabic{theorem}}
\renewcommand{\theHequation}{appendix.\arabic{equation}}
\renewcommand{\thesubsection}{\arabic{subsection}}
\renewcommand{\thetheorem}{\arabic{theorem}}
\renewcommand{\theequation}{\arabic{equation}}
\section*{Appendix. Sharpness and counterexamples}
\label{sec:examples}\label{sec:obstructions}
We give examples showing the sharpness and limitations of the
preceding estimates. For the first example, we verify $k$-convexity
and identify the full Hessian measure in detail. These arguments
will also be used in the subsequent constructions.

\setcounter{subsection}{0}
\setcounter{theorem}{0}
\setcounter{equation}{0}
\subsection{Critical densities and logarithmic moduli}
\begin{proposition}\label{prop:criticaldensity}
For every $2\le k<n$, there is a function
\[
 u\in\Phi_k(B_1)\cap C(B_1)\cap W^{2,k}_{\mathrm{loc}}(B_1),
 \quad F_k[u]=f\,dx,
\]
smooth and strictly $k$-convex away from the origin, such that
\[
 f\in L^{n/k}(B_1),\quad
 |Du(x)|\to\infty\text{ as }x\to0.
\]
\end{proposition}
\begin{proof}
\emph{Step 1. }
Set $L(r)=\log(4/r)$, $g(r)=\log L(r)$,
$H(r)=\int_0^rL(t)^{-1}\,dt$, and $A_0=6n/(n-k)$, and set
\begin{equation}\label{eq:criticalcounterexample}
 u(x)=x_1g(|x|)+A_0H(|x|)\quad(x\ne0),\quad u(0)=0.
\end{equation}
Since $r\log L(r)\to0$, $u$ is continuous. For $e=x/r$ and
$\beta=(rL(r))^{-1}$,
\[
 D^2H=\beta[I+(L^{-1}-1)e\otimes e],\quad
 \|D^2(x_1g)\|_{\mathrm{op}}\le3\beta.
\]
In a radial frame, $D^2u\ge\beta M_0$, where
$M_0=\operatorname{diag}(-3,A_0-3,\ldots,A_0-3)$. For $1\le j\le k$,
\[
 \sigma_j(M_0)=(A_0-3)^{j-1}
 \left[\binom{n-1}{j}(A_0-3)-3\binom{n-1}{j-1}\right]>0,
\]
because $A_0>3n/(n-k)\ge3n/(n-j)$. Cone monotonicity and
$|D^2u|\le C\beta$ therefore yield
\begin{equation}\label{eq:criticalcounterdensity}
 D^2u\in\Gamma_k,\quad
 c(rL(r))^{-k}\le f:=\sigma_k(D^2u)\le C(rL(r))^{-k}.
\end{equation}

\emph{Step 2. }
Here $|Du|\le g+A_0+1$ and
$\int_0^1r^{n-k-1}L(r)^{-k}\,dr<\infty$. The inner boundary terms
in integration by parts for the first and second derivatives are
$O(r^ng(r))$ and $O(r^{n-1}(g(r)+1))$, respectively. They vanish,
so $u\in W^{2,k}_{\mathrm{loc}}$ with its classical Hessian.

We record the measure argument used also below. If continuous
$v\in W^{2,k}_{\mathrm{loc}}$ satisfies
$D^2v\in\overline{\Gamma_k}$ almost everywhere, local convolution
preserves this closed convex cone and converges locally uniformly to
$v$. Hence $v\in\Phi_k$. Moreover,
\[
 |\sigma_k(A)-\sigma_k(B)|
 \le C(|A|+|B|)^{k-1}|A-B|
\]
and strong $L^k$ convergence of the Hessians imply strong $L^1$
convergence of their densities. Weak continuity of Hessian measures
then gives $F_k[v]=\sigma_k(D^2v)\,dx$.
For \eqref{eq:criticalcounterexample}, this excludes any singular mass.

\emph{Step 3. }
By \eqref{eq:criticalcounterdensity},
\[
 \int_{B_1}f^{n/k}\,dx\le C\int_0^1\frac{dr}{rL(r)^n}<\infty\text{ and }
 Du=e_1\log L(r)+\left(A_0-\frac{x_1}{r}\right)\frac{x}{rL(r)}.
\]
The second gradient term tends to zero, while the first diverges.
\end{proof}

\begin{proposition}\label{obs:loglip-sharp}
For every $2\le k<n$, there is a continuous weak $k$-convex function $u$ on $B_1$, smooth and strictly $k$-convex away from zero, such that $F_k[u]=f\,dx$ with
\[
 0<f(x)\asymp |x|^{-k},\quad
 \sup_{z\in B_{1/2},\ 0<\rho<1/4}
       \frac{F_k[u](B_\rho(z))}{\rho^{n-k}}<\infty.
\]
The function has a local modulus $C\rho\log(4/\rho)$, and this order cannot be replaced by $o(\rho\log(1/\rho))$ at zero.
\end{proposition}

\begin{proof}
Set $u(x)=x_1\log(4/|x|)+A_0|x|$ for $x\ne0$, $u(0)=0$,
and $A_0=6n/(n-k)$. With $e=x/r$ and $\xi=x_1/r$,
\[
 D^2(x_1\log(4/r))
 =r^{-1}[-e_1\otimes e-e\otimes e_1+\xi(2e\otimes e-I)],\quad
 D^2(A_0r)=A_0r^{-1}(I-e\otimes e).
\]
The first matrix has norm at most $3/r$. The comparison in
Proposition~\ref{prop:criticaldensity}, Step 1, gives
$D^2u\in\Gamma_k$ and $f=\sigma_k(D^2u)\asymp r^{-k}$.

Since $|D^2u|\le C/r$ and $|Du|\le C\log(4/r)$, Step 2 of that
proof applies: $u\in W^{2,k}_{\mathrm{loc}}\cap\Phi_k$ and
$F_k[u]=f\,dx$, without an atom at zero. The bound
\[
 \int_{B_\rho(z)}|x|^{-k}\,dx\le C\rho^{n-k}
\]
is uniform in $z$: compare with $B_{3\rho}(0)$ for $|z|\le2\rho$
and use $|x|\ge|z|/2$ otherwise.

The upper modulus follows from
Corollary~\ref{pot:log-lipschitz}; its order is attained because
$u(te_1)-u(0)=t\log(4/t)+A_0t$ for $t>0$.
\end{proof}

\subsection{Failure of gradient BMO bounds}\label{sharp:bmo-failure}\label{obs:morrey-subsection}\label{obs:critical-subsection}
For a locally integrable vector field $F$, set
\[
 E_F(x,r)=\avint_{B_r(x)}|F-(F)_{B_r(x)}|.
\]
The two constructions below separate critical size conditions from
locally bounded mean oscillation of the gradient.

\begin{proposition}\label{obs:bmo-limits}
\begin{enumerate}
\renewcommand{\labelenumi}{\textup{(\roman{enumi})}}
\renewcommand{\theenumi}{\thetheorem(\roman{enumi})}
\item\label{obs:folded-morrey}
For every $2\le k<n$, there is
$v\in\Phi_k(B_1)\cap C(B_1)$ with
$Dv\notin\mathrm{BMO}_{\mathrm{loc}}(B_1)$ and
\begin{equation}\label{obs:uniform-vanishing-morrey}
 \lim_{\rho\to0}\sup_{z\in B_{1/2}}
 \frac{F_k[v](B_\rho(z))}{\rho^{n-k}}=0.
\end{equation}
Its measure has a bulk density and a part on $\{x_1=0\}$, with no
atom at the origin.
\item\label{obs:adaptive-critical}
For every $3\le k<n$, there is
$u\in\Phi_k(B_1)\cap C(B_1)\cap W^{2,k}_{\mathrm{loc}}(B_1)$
with $F_k[u]=f\,dx$, $f\in L^{n/k}(B_1)$, and
$Du\notin\mathrm{BMO}_{\mathrm{loc}}(B_1)$.
The function is smooth and strictly $k$-convex away from the origin.
\end{enumerate}
\end{proposition}
\begin{proof}
\emph{Part (i). }
Use the notation
\begin{equation}\label{obs:folded-definitions}
 L(r)=\log(4/r),\quad a(r)=\log L(r),\quad b(r)=L(r)^{-1},\quad
 H(r)=\int_0^rb(t)\,dt,\quad A_0=\frac{6n}{n-k}.
\end{equation}
Let $g=x_1a(|x|)$, $h=A_0H(|x|)$, and $v=h+|g|$, with $v(0)=0$.
Proposition~\ref{prop:criticaldensity} gives $h\pm g\in\Phi_k$;
their maximum $v$ is therefore weakly $k$-convex and continuous.

On the half-balls, the Sobolev chain rule gives
\begin{equation*}
 Dv=\operatorname{sgn}(x_1)e_1a(r)
       +\left(A_0-\frac{|x_1|}{r}\right)b(r)e,\quad e=x/r.
\end{equation*}
This gradient is locally integrable. Since $v$ is even,
$(Dv)_{B_r}=0$; projection onto $\operatorname{sgn}(x_1)e_1$ yields
\begin{equation*}
 E_{Dv}(0,r)=\avint_{B_r}|Dv|\ge a(r)=\log L(r)\to\infty.
\end{equation*}

Writing $x=(x_1,y)$, we claim
\begin{equation}\label{obs:folded-measure}
 \begin{split}
 F_k[v]&=f\,dx+w(|y|)\,\mathcal H^{n-1}\restr\{x_1=0\},\quad
 f(x)=\sigma_k(D^2(h+\operatorname{sgn}(x_1)g))(x)\quad(x_1\ne0),\\
 w(r)&=\frac{2A_0^{k-1}a(r)}{(rL(r))^{k-1}}
 \left[\binom{n-2}{k-1}+\frac1{L(r)}\binom{n-2}{k-2}\right].
 \end{split}
\end{equation}
In particular,
\begin{equation*}
 0\le f(x)\le C(rL(r))^{-k},\quad
 w(r)\asymp\frac{\log L(r)}{(rL(r))^{k-1}}.
\end{equation*}

To identify the interface term, set
$\phi_\delta(t)=\sqrt{t^2+\delta^2}-\delta$ and
$v_\delta=h+\phi_\delta(g)$. Then
\begin{equation*}
 D^2v_\delta=B_\delta+\phi_\delta''(g)Dg\otimes Dg,
 \quad B_\delta=D^2h+\phi_\delta'(g)D^2g,
\end{equation*}
where $B_\delta\in\Gamma_k$ by convexity of the cone, and
\begin{equation*}
 |B_\delta|\le Cb/r,\quad |Dg|\le Ca.
\end{equation*}
For fixed $\delta$, these bounds and
$\phi_\delta''(g)|Dg|^2\le\delta^{-1}|Dg|^2$ imply
$v_\delta\in W^{2,k}_{\mathrm{loc}}$; the logarithmic growth of its
gradient makes the inner integration-by-parts terms vanish.
Proposition~\ref{prop:criticaldensity}, Step 2, therefore identifies
its full measure with its density. Also $|v_\delta-v|\le\delta$, so
these measures converge weakly to $F_k[v]$.

The Newton tensor rank-one identity is
\begin{equation}\label{obs:constant-fold-measure}
 \sigma_k(D^2v_\delta)=\sigma_k(B_\delta)
 +\phi_\delta''(g)\langle T_{k-1}(B_\delta)Dg,Dg\rangle.
\end{equation}

On compact portions of $\{x_1=0\}\setminus\{0\}$,
$Dg=ae_1$, $D_y^2g=0$, and
$(T_{k-1}(B_\delta))_{11}=\sigma_{k-1}(D_y^2h)$.
Work in a level-set coordinate neighborhood of the interface,
away from the origin, and let $\psi$ be supported there.
The coarea formula writes the second term as
\[
 \int_{\mathbb R}\phi_\delta''(t)
 \int_{\{g=t\}}\psi\,
 \frac{\langle T_{k-1}(B_\delta)Dg,Dg\rangle}{|Dg|}
 \,d\mathcal H^{n-1}\,dt.
\]
As the level parameter $t\to0$, the inner integral converges,
uniformly in $\delta$, to the integral of
$\psi a\sigma_{k-1}(D_y^2h)$ over the interface.
Indeed, $Dg\to ae_1$, the tangential block $D_y^2g\to0$, and
$|\phi_\delta'|\le1$; only that tangential block enters
$(T_{k-1})_{11}$.

Since $\phi_\delta''(t)\,dt\to2\delta_0$ weakly as Radon measures, the limiting
surface density is $2a\sigma_{k-1}(D_y^2h)$.
The tangential eigenvalues $A_0/(rL^2)$ and $A_0/(rL)$ (multiplicity
$n-2$) give $w$ above. Dominated convergence gives the bulk term
on compact sets away from zero.

The bulk and interface measures have now been identified away
from zero, so the only remaining possible singular part is an atom
at the origin. To exclude it, use \eqref{obs:constant-fold-measure}:

\[
 \sigma_k(D^2v_\delta)\le C(b/r)^k+
 C(b/r)^{k-1}a^2\phi_\delta''(x_1a).
\]
Since $n\ge3$,
\[
 \int_{S^{n-1}}\phi_\delta''(ra\theta_1)\,d\theta
 \le C\int_{-1}^1\phi_\delta''(ra t)\,dt\le C/(ra).
\]
Consequently, uniformly in $\delta$,
\begin{equation*}
 F_k[v_\delta](B_\rho)
 \le C\int_0^\rho s^{n-k-1}[b(s)^k+a(s)b(s)^{k-1}]\,ds
 \le C\rho^{n-k}\frac{\log L(\rho)}{L(\rho)^{k-1}}.
\end{equation*}
The last weights are increasing for sufficiently small $s$ and
$n>k$. The Portmanteau inequality on $B_\rho$ gives the same bound
for $F_k[v]$; it tends to zero, so \eqref{obs:folded-measure} is
the full measure.

Its surface term on $\rho/2<|y|<\rho$ also gives
\[
 F_k[v](B_\rho)\asymp\rho^{n-k}\log L(\rho)/L(\rho)^{k-1}.
\]

The estimates
\[
 \int_{B_\rho(z)}|x|^{-k}\,dx\le C\rho^{n-k},\quad
 \int_{B_\rho(z)\cap\{x_1=0\}}|y|^{-(k-1)}\,d\mathcal H^{n-1}
 \le C\rho^{n-k}
\]
hold for all centers, by the near/far split used in
Proposition~\ref{obs:loglip-sharp}, also in dimension $n-1$.

Given $\eta>0$, both weights $b^k$ and $ab^{k-1}$ are at most
$\eta$ in a sufficiently small fixed ball. Its contribution to
$F_k[v](B_\rho(z))/\rho^{n-k}$ is at most $C\eta$.
Outside it the volume and surface densities are bounded, contributing
$C_\eta\rho^k$ and $C_\eta\rho^{k-1}$. Letting $\rho\to0$
and then $\eta\to0$ proves \eqref{obs:uniform-vanishing-morrey}.

\medskip\noindent\emph{Part (ii). }
Keep \eqref{obs:folded-definitions}, $g=x_1a(|x|)$, and $h=A_0H(|x|)$.
The variable smoothing width
\begin{equation*}
 \varepsilon(r)=rb(r),\quad F=(g^2+\varepsilon^2)^{1/2},\quad
 u=h+F,\quad u(0)=0
\end{equation*}
removes the interface measure. Since $F\le r(a+b)$, $u$ is continuous
and smooth away from zero.

Also $\varepsilon'=b+b^2$ and
$\varepsilon''=(b^2+2b^3)/r$, so $D^2\varepsilon\ge0$ and
$|D^2\varepsilon|\le Cb/r$. Set $\tau=g/F$ and $\eta=\varepsilon/F$;
then $|\tau|\le1$, $0\le\eta\le1$, and $\tau^2+\eta^2=1$.

The Hessian splits into an admissible base matrix and a nonnegative
rank-one term:
\begin{equation*}
 D^2u=B+Q,\quad B=D^2h+\tau D^2g+\eta D^2\varepsilon,\quad
 Q=\frac{(\varepsilon Dg-gD\varepsilon)^{\otimes2}}{F^3}.
\end{equation*}
Here $Q\ge0$, $\eta\ge0$, and the first two terms of $B$ are a
convex combination of $D^2(h\pm g)$. Hence
\begin{equation*}
 B,B+Q\in\Gamma_k,\quad |B|\le Cb/r.
\end{equation*}

For $e=x/r$ and $\xi=x_1/r$, $Dg=ae_1-\xi be$ and
$D\varepsilon=(b+b^2)e$, whence
\begin{equation*}
 \varepsilon Dg-gD\varepsilon
 =rab\left[e_1-\xi\left(1+b+\frac ba\right)e\right].
\end{equation*}
The bracket is bounded, and $F=rb(1+a^2\xi^2/b^2)^{1/2}$, so
\begin{equation*}
 |Q|\le\frac{Ca^2}{rb}\left(1+\frac{a^2\xi^2}{b^2}\right)^{-3/2}.
\end{equation*}
The largest second derivatives lie in the angular strip
$|x_1|/r\lesssim b/a$. Its width provides the gain in the next integral.

For $p=n/k>1$, the spherical density in $\xi$ is bounded ($n\ge3$).
Integration over $\xi$, followed by $s=a\xi/b$, gives
\begin{equation*}
 \int_{S^{n-1}}|Q|^p\,d\theta
 \le C\left(\frac{a^2}{rb}\right)^p
 \int_{-1}^1(1+a^2\xi^2/b^2)^{-3p/2}\,d\xi
 \le C\left(\frac{a^2}{rb}\right)^p\frac ba.
\end{equation*}
Since $Q$ has rank one, it enters $\sigma_k(B+Q)$ only linearly.
The rank-one identity therefore bounds the classical density by
\[
 0<f=\sigma_k(D^2u)\le C(b/r)^k+C(b/r)^{k-1}|Q|.
\]
Since $kp=n$, it follows that
\begin{equation*}
 \int_{B_1}f^p\,dx
 \le C\int_0^1[b^{kp}+a^{2p-1}b^{(k-2)p+1}]\,\frac{dr}{r}
 =C\int_{\log4}^\infty
 [L^{-n}+(\log L)^{2p-1}L^{-((k-2)p+1)}]\,dL<\infty.
\end{equation*}
The last inequality uses $(k-2)p+1>1$, precisely where $k\ge3$
is needed.

The coarse bound $|D^2u|\le Cr^{-1}L(r)(\log L(r))^2$ is in
$L^k$ because $n>k$. Also $|Du|\le a+Cb$ and $|u|\le Cr(a+b)$.

The inner boundary terms for the first and second weak derivatives
are respectively $O(r^n(a+b))$ and $O(r^{n-1}(a+b))$, and vanish.
Thus $u\in W^{2,k}_{\mathrm{loc}}$. The measure argument in
Proposition~\ref{prop:criticaldensity}, Step 2, now gives
$u\in\Phi_k$ and $F_k[u]=f\,dx$, including at the origin.

Finally, $u$ is even, $(Du)_{B_r}=0$, and
\[
 Du=\tau ae_1+[A_0b-\tau\xi b+\eta(b+b^2)]e.
\]
The bracket is at least $(A_0-1)b>0$. On the fixed angular set
$|\xi|\ge1/2$, $|\tau|\ge1/2$ at small radii since $a/b\to\infty$.
Projection onto $\operatorname{sgn}(\xi)e_1$ therefore gives
$|Du|\ge a/2$. This set has a fixed positive relative area and
$a(s)\ge a(r)$ for $s<r$, so
\[
 E_{Du}(0,r)\ge c\log L(r)\to\infty.
\]

\end{proof}

\subsection{Limitations of the potential estimate}\label{sharp:critical-potential}
\begin{proposition}
\label{rem:radialpotential}\label{obs:riesz-root}
For every $2\le k<n$:
\begin{enumerate}
 \item[(i)] There is a convex $v\in C^1(B_1)$ with bounded gradient
 whose potential $\W_k^R[F_k[v]]$ is unbounded near zero, for every
 fixed $0<R<1/2$.
 \item[(ii)] The Wolff term in \eqref{eq:main} cannot be replaced by
 $(I_k^R[F_k[u]])^{1/k}$, where
 \begin{equation*}
  I_k^R[\mu](x)=\int_0^R\frac{\mu(B_t(x))}{t^{n-k}}\,\frac{dt}{t}.
 \end{equation*}
\end{enumerate}
\end{proposition}
\begin{proof}
\emph{(i)} Use $L(r)=\log(4/r)$ and $H(r)=\int_0^rL(t)^{-1}\,dt$
from Proposition~\ref{prop:criticaldensity}, and set $v(x)=H(|x|)$.
Then $|Dv|=1/L(r)\to0$, and the radial and tangential Hessian
eigenvalues are $(rL^2)^{-1}$ and $(rL)^{-1}$. Also $H(r)/r\le1/L(r)\to0$, so $v$ is convex
and $C^1$, with $Dv(0)=0$.

The full-measure identification in the
proof of Proposition~\ref{prop:criticaldensity} applies, giving
$F_k[v]=f_v\,dx$ with $f_v\asymp(rL)^{-k}$ and no mass at zero.
Hence
\[
 F_k[v](B_t)\asymp t^{n-k}L(t)^{-k}.
\]

In particular the
potential is infinite at zero; for $0<|x|<R/4$, the inclusion
$B_{t/2}(0)\subset B_t(x)$ for $t\ge2|x|$ gives
\[
 \W_k^R[F_k[v]](x)\ge c\int_{2|x|}^R\frac{dt}{tL(t)}
 \to\infty\quad(x\to0).
\]

\emph{(ii)} Take $u$ from Proposition~\ref{prop:criticaldensity}.
Its full measure is $f\,dx$ with $f(y)\le C(|y|L(|y|))^{-k}$.
Tonelli's theorem gives
\[
 I_k^{1/2}[f\,dx](x)
 \le C\int_{B_1} f(y)|x-y|^{k-n}\,dy.
\]
For $0<d=|x|<1/4$, split the last integral into
$|y|<d/2$, $d/2\le|y|\le2d$, and $|y|>2d$. The first two
contributions are each at most $CL(d)^{-k}$: use respectively
$|x-y|\ge d/2$ and
$f(y)\le Cd^{-k}L(d)^{-k}$ together with
$\int_{B_{3d}(x)}|x-y|^{k-n}\,dy\le Cd^k$.

The third contribution is bounded by
\[
 C\int_{2d}^1dt/(tL(t)^k).
\]
This integral remains finite with lower
limit zero, because $k>1$, also covering $x=0$. Therefore
\begin{equation*}
 \sup_{x\in B_{1/4}}I_k^{1/2}[F_k[u]](x)<\infty.
\end{equation*}
For the fixed radius $R=1/2$, the average term in \eqref{eq:main}
is uniformly bounded as the center tends to zero, whereas
$|Du(x)|\to\infty$ through nonzero gradient Lebesgue points.
This contradicts the proposed replacement.
\end{proof}

\subsection{Failure of uniform boundary gradient bounds}\label{sharp:boundary}
The following sequence shows why an interior estimate cannot be
extended uniformly to the boundary using only the same measure data
and the height of the solution.

\begin{proof}
Write \(a=t+\varepsilon\), suppress subscripts, and set
\(F=L^{1/(4k)}\). Since \(L>3\),
\[
F'=-\frac{F}{4kaL},\quad
F''=\frac{|F'|}{a}\left(1-\frac{1-1/(4k)}{L}\right),
\quad
\frac{2(F')^2}{F}<F''\le\frac{|F'|}{a}.
\]
In the orthogonal coordinates \((y,t)\), the Hessian has blocks
\[
D^2w_\varepsilon=
\begin{pmatrix}
2FI_{n-1}&2F'y\\
2F'y^T&c
\end{pmatrix},
\quad
c=2F+4(t-1)F'+(|y|^2-2t+t^2)F''.
\]
Its Schur complement satisfies
\[
\begin{aligned}
S
&=2F+4(t-1)F'-t(2-t)F''
+|y|^2\left(F''-\frac{2(F')^2}{F}\right)\\
&\ge 2F+\frac{|F'|}{a}\bigl[4(1-t)a-t(2-t)\bigr]
\ge 2F-2a|F'|
=2F\left(1-\frac{1}{4kL}\right)>0.
\end{aligned}
\]
Here
\[
4(1-t)a-t(2-t)+2a^2=t(2-t)+4\varepsilon+2\varepsilon^2\ge0.
\]
Hence \(w_\varepsilon\) is strictly convex.

Expanding the principal minors, the off-diagonal contribution is
nonpositive. Since \(c\le2F+4|F'|\), \(F^k=L^{1/4}\), and
\(F^{k-1}|F'|=1/(4kaL^{3/4})\), we obtain
\[
0<f_\varepsilon^0:=\sigma_k(D^2w_\varepsilon)
\le C(F^k+F^{k-1}|F'|)
\le\frac{C}{aL^{3/4}}
\le\frac{C}{tL_0(t)^{3/4}}.
\]
Here \(aL\) is bounded, \(aL^{3/4}\) is increasing on \((0,9/4]\),
and \(L_0(t)=\log(e^4/t)\). The horizontal cross-section of \(B_1\)
has measure \(\omega_{n-1}(2t-t^2)^{(n-1)/2}\). Consequently, writing
\(q=q_\partial\),
\[
\sup_\varepsilon\int_{B_1}(f_\varepsilon^0)^q\,dx
\le C\left(1+\int_0^1\frac{dt}{tL_0(t)^{3q/4}}\right)
<\infty,
\]
because \(q\ge2\).

Set \(c_\varepsilon=F(0)^{-1/2}\),
\(u_\varepsilon=c_\varepsilon w_\varepsilon\), and
\(f_\varepsilon=c_\varepsilon^k f_\varepsilon^0\).
The bound \(|w_\varepsilon|\le2tL_0(t)^{1/(4k)}\) is uniform, so
\(c_\varepsilon\to0\) gives vanishing height and density norms, while
\[
 \partial_\nu u_\varepsilon(e_n)=2F(0)^{1/2}\to\infty.
\]

Since \(q>n/k\), Hölder's inequality gives
\[
\sup_{x\in\mathbb R^n}\W_k^r[f_\varepsilon\mathbf1_{B_1}\,dx](x)
\le C\|f_\varepsilon\|_{L^q(B_1)}^{1/k}
r^{1-n/(kq)},\quad 0<r\le4.
\]
Taking \(r=4\) proves the vanishing potential norms. The positive power
of \(r\), together with the uniform density bound, also gives uniformly
vanishing small-scale tails. Letting \(\varepsilon\to0\) along any
sequence proves the claim.
\end{proof}
\noindent \textbf{Declaration of AI-assisted research and writing.} We used ChatGPT and Codex (OpenAI) to assist with mathematical
exploration and with drafting and revising this manuscript.
In particular, these tools assisted with some counterexample
constructions and with the design of the two-point test function
used in the proof of Proposition~\ref{prop:defect}. The authors take full responsibility for the mathematical
content and the final manuscript.

\end{document}